\documentclass[final,onefignum,onetabnum]{siamart220329}

\usepackage{lipsum}
\usepackage{amsfonts}
\usepackage{graphicx}
\usepackage{epstopdf}
\usepackage{algorithmic}
\ifpdf
  \DeclareGraphicsExtensions{.eps,.pdf,.png,.jpg}
\else
  \DeclareGraphicsExtensions{.eps}
\fi

\usepackage{enumitem}

\newsiamremark{remark}{Remark}

\headers{Domain decomposition-based training algorithms for PINNs}{A. Kopani\v{c}\'akov\'a, H. Kothari, G. Em Karniadakis, and R. Krause}

\title{Enhancing Training of Physics-Informed Neural Networks using Domain-Decomposition based Preconditioning Strategies\thanks{%Submitted to the editors 30-June-2023.
\funding{
A.K. was supported by the Swiss National Science Foundation (SNSF) under the project ``\textit{Multilevel training of DeepONets -- multiscale and multiphysics applications}''.
A.K., H.K., and R.K. acknowledge the support of the Platform for Advanced Scientific Computing (PASC) under the project \textit{ExaTrain} as well as the support of the SNSF under the project ``\textit{ML$^2$ - Multilevel and Domain Decomposition Methods for Machine Learning}''. G.E.K. would like to acknowledge the support of the Vannevar Bush Faculty Fellowship.}}}

\author{
Alena Kopani\v{c}\'akov\'a\footnotemark[3]
\thanks{Division of Applied Mathematics, Brown University, Providence, USA
   (\email{alena\_kopanicakova@brown.edu, george\_karniadakis@brown.edu}).}
\and Hardik Kothari\thanks{Euler Institute, Universit\`a della Svizzera italiana, Lugano, Switzerland
  (\email{alena.kopanicakova@usi.ch}, \email{hardik.kothari@usi.ch}, \email{rolf.krause@usi.ch}).}
\and George Em Karniadakis\footnotemark[2]
\and Rolf Krause\footnotemark[3]
\thanks{UniDistance, Switzerland
   (\email{rolf.krause@fernuni.ch}).}
}

\usepackage{amsopn}

\ifpdf
\hypersetup{
  pdftitle={Domain decomposition-based training algorithms for PINNs},
  pdfauthor={Alena Kopani\v{c}\'{a}kov\'{a}, Hardik Kothari, George Em Karniadakis and Rolf Krause}
}
\fi

\usepackage[T1]{fontenc}
\usepackage[latin1]{inputenc}
\usepackage{caption}
\usepackage{algorithm}
\usepackage{graphicx,url}
\usepackage{amsmath,amssymb,amsopn}
\usepackage{tabularx}
\usepackage{color}
\usepackage{calc}
\usepackage{longtable}
\usepackage{multirow}
\usepackage{mathrsfs}
\usepackage{array}
\usepackage{hyperref}
\usepackage{tikz}
\usepackage{booktabs}
\usepackage{pgfplotstable}
\usepackage{pgfplots}
\usepackage[outline]{contour}
\usepackage[colorinlistoftodos,prependcaption]{todonotes}
\usepackage{calrsfs}
\usepackage{graphicx}
\usepgfplotslibrary{groupplots}
\usetikzlibrary{matrix}
\usetikzlibrary{automata, external, shapes.geometric, fit}
\usetikzlibrary{arrows.meta}
\pgfplotsset{compat=1.13}

\def \bv{\vec{b}}

\def \gv{\vec{g}}

\def \pv{\vec{p}}

\def \sv{\vec{s}}

\def \vv{\vec{v}}

\def \xv{\vec{x}}
\def \yv{\vec{y}}

\def \rm{\mat{r}}

\def \Dm{\mat{D}}

\def \Im{\mat{I}}

\def \Lm{\mat{L}}

\def \R{\mathbb{R}}			\def \N {\mathbb{N}}

\DeclareMathAlphabet{\pazocal}{OMS}{zplm}{m}{n}

\renewcommand{\vec}[1]{\boldsymbol{#1}}
\newcommand{\mat}[1]{\boldsymbol{\mathrm{#1}}}

\definecolor{myblack}{RGB}{53, 53, 53}
\definecolor{myblue}{RGB}{40, 75, 99}
\definecolor{myred}{RGB}{192, 50, 33}
\definecolor{myyellow}{RGB}{255, 166, 48}
\definecolor{mywhite}{RGB}{240, 237, 238}
\definecolor{mygreen}{RGB}{0, 102, 0}

\definecolor{green1}{RGB}{9, 82, 86}
\definecolor{green2}{RGB}{8, 127, 140}
\definecolor{green3}{RGB}{6, 167, 125}
\definecolor{green4}{RGB}{79, 109, 122}
\definecolor{green5}{RGB}{192, 214, 223}
\definecolor{violet}{RGB}{26,69,131}

\definecolor{checkgreen}{rgb}{0,0.6,0}
\definecolor{phase1}{rgb}{0.008,0.655,1.000}
\definecolor{phase2}{rgb}{0.016,0.75,0.700}
\definecolor{phase3}{rgb}{0.929,0.35,0.700}
\definecolor{icsyellow}{cmyk}{0.00,0.11,0.53,0.00}

\usepackage{color, colortbl}
\definecolor{Gray}{gray}{0.9}

\definecolor{checkgreen}{rgb}{0,0.6,0}
\definecolor{phase1}{rgb}{0.008,0.655,1.000}
\definecolor{phase2}{rgb}{0.016,0.75,0.700}
\definecolor{phase3}{rgb}{0.929,0.35,0.700}
\definecolor{icsyellow}{cmyk}{0.00,0.11,0.53,0.00}
\definecolor{green1}{RGB}{9, 82, 86}
\definecolor{green2}{RGB}{8, 127, 140}
\definecolor{green3}{RGB}{6, 167, 125}
\definecolor{green4}{RGB}{79, 109, 122}
\definecolor{green5}{RGB}{192, 214, 223}
\definecolor{violet}{RGB}{26,69,131}
\definecolor{mypurple}{RGB}{150, 0, 180}

\usetikzlibrary{intersections}
\usepgfplotslibrary{fillbetween}
\usepackage{multirow}
\usepackage{stmaryrd}
\usepackage{xfrac}
\usepackage{caption}
\usepackage{subcaption}
\usepackage{lmodern}
\usetikzlibrary{decorations.text}
\usepackage[normalem]{ulem}

\definecolor{blackmy}{RGB}{38, 70, 83}
\definecolor{greenmy}{RGB}{42, 167, 143}
\definecolor{yellowmy}{RGB}{233, 196, 106}
\definecolor{redmy}{RGB}{253,127,111}

\definecolor{bluemy}{RGB}{40, 75, 99}
\definecolor{myred}{RGB}{192, 50, 33}
\definecolor{brownmy}{RGB}{244, 162, 97}

\definecolor{blue1}{RGB}{1, 58, 99}
\definecolor{blue2}{RGB}{70, 143, 175}
\definecolor{blue3}{RGB}{137, 194, 217}

\definecolor{red1}{RGB}{176, 71, 89}
\definecolor{red3}{RGB}{249, 155, 125}
\definecolor{red2}{RGB}{231, 97, 97}

\definecolor{green1}{RGB}{112, 151, 117}
\definecolor{green2}{RGB}{143, 185, 150}
\definecolor{green3}{RGB}{161, 204, 165}

\DeclareMathOperator*{\argmin}{arg\,min}

\newcommand{\algorithmiccommentMine}[1]{\bgroup\hfill$\triangleright$~{#1}\egroup}
\newcommand\COMMENTmine[1]{\algorithmiccommentMine{#1}}

\newcommand\newtext[1]{#1}
\newcommand\newtextreviewerone[1]{#1}
\newcommand\oldtext[1]{}
\newcommand\cancel[1]{}

\begin{document}

\maketitle

\begin{abstract}
	We propose to enhance the training of physics-informed neural networks (PINNs).
	To this aim, we introduce nonlinear additive and multiplicative preconditioning strategies for the widely used L-BFGS optimizer.
	The nonlinear preconditioners are constructed by utilizing the Schwarz domain-decomposition framework, where the parameters of the network are decomposed in a layer-wise manner.
	Through a series of numerical experiments, we demonstrate that both, additive and multiplicative preconditioners significantly improve the convergence of the standard L-BFGS optimizer, while providing more accurate solutions of the underlying partial differential equations.
	Moreover, the additive preconditioner is inherently parallel, thus giving rise to a novel approach to model parallelism.
\end{abstract}



\section{Introduction}
In recent years, deep learning approaches for solving partial differential equations (PDEs) have gained popularity due to their simplicity, mesh-free nature, and efficiency in solving inverse and high-dimensional problems~\cite{karniadakis2021physics, weinan2021dawning}.
The notable deep learning approaches in the context of scientific machine learning include physics-informed neural networks (PINNs)~\cite{raissi2019physics}, deep Galerkin method~\cite{sirignano2018dgm}, and deep Ritz method~\cite{yu2018deep}.
The idea behind these methods is to utilize deep neural networks (DNNs) to approximate solutions of the PDEs by minimizing a loss function, which incorporates the physical laws and domain knowledge.
More specifically, in the context of PINNs, which we consider in this work, the loss function incorporates a PDE residual, boundary/initial conditions, and sometimes also observational data.

To approximate a solution of a PDE using PINNs, the optimal parameters of DNN have to be found during the training.
This is typically carried out using Adam~\cite{kingma2014adam} and/or L-BFGS~\cite{liu1989limited} optimizers.
However, it has been shown that the performance of these optimizers deteriorates for highly ill-conditioned and stiff problems~\cite{wang2021understanding}.
Several algorithmic enhancements have been therefore proposed to alleviate such difficulties, for instance, adaptive weighing strategies~\cite{wang2022and}, empirical learning-rate annealing scheme~\cite{wang2021understanding}, multiscale architectures~\cite{dolean2023multilevel, gratton2023block}, ensemble learning~\cite{fang2023ensemble} based on gradient boosting~\cite{callens2020using}, adaptive activation functions~\cite{jagtap2020adaptive, jagtap2020locally} etc.
Furthermore, novel optimizers have also been developed, for instance, Gauss-Newton methods~\cite{schraudolph2002fast, botev2017practical}, Levenberg--Marquardt optimizers~\cite{more1978levenberg, yu2018levenberg, gavin2019levenberg}, variable projection methods~\cite{newman2020train, golub2003separable, o2013variable}, block-coordinate methods~\cite{cyr2020robust}, etc.

In this work, our goal is to enhance the convergence rate of the optimizers used for the training of PINNs by developing novel nonlinear preconditioning strategies.
In the field of scientific computing, such nonlinear preconditioning strategies are often used to accelerate the convergence of the nonlinear iterative methods~\cite{brune2015composing}.
Specifically, left preconditioners~\cite{cai2002nonlinearly, liu2015field} enhance the conditioning of nonlinear problems by rebalancing nonlinearities or by transforming the basis of the solution space.
In contrast, right preconditioners~\cite{luo2023preconditioned,klawonn2023adaptive} provide an improved initial guess for a given iterate.
Both, left and right, preconditioners can be constructed using various approaches.
Prominent approaches typically rely on the decomposition of the solution space into subspaces, such as field-split~\cite{kopanivcakova2023nonlinear,liu2015field}, domain-decomposition (DD)~\cite{cai2002nonlinearly, kothari2022nonlinear,kothari2022nonlinear_bounds}, or multilevel decomposition~\cite{yuan2014review,kopanicakova2022globally, Kopanicakova_2020c,gratton2023multilevel,kopanivcakova2019recursive}.
Among these approaches, the DD approach is of particular interest due to its parallelization capabilities.
As a consequence, preconditioners constructed using the DD approach not only improve the convergence rate of a given iterative method but also allow for a reduction of the overall computation time, as the workload can be distributed between larger amounts of computational resources.

In the field of machine learning, the parallelization of the training algorithms has attracted a lot of attention recently, see~\cite{dean2012large, mnih2016asynchronous, ben2019demystifying} for a comprehensive overview.
The developed approaches can be broadly classified into two main groups: data parallel and model parallel.
The data parallel approach has gained significant popularity due to its architecture-agnostic nature.
In contrast, model parallelism is more intricate as it involves decomposing the network across different compute nodes.
As a consequence, this approach has to be typically tailored to a specific network architecture~\cite{ben2019demystifying}.
For instance, the {\sf TorchBraid} framework~\cite{gunther2020layer} enables layer-parallel training of ODE-based networks by utilizing a multigrid-in-time approach.
In~\cite{gu2023decomposition, gu2022decomposition}, a model parallel approach for convolutional networks is proposed, where the network's width is partitioned into subnetworks, which can be trained in parallel.
The parameters obtained from training these subnetworks are then used to initialize the global network, which is further trained sequentially.
Similarly, in~\cite{klawonn2023domain}, a DD-based network architecture is proposed for parallel training of image recognition tasks.

In the context of PINNs, DD-based approaches have also been utilized to enhance training efficiency and scalability to larger computational domains and multi-scale problems.
These approaches leverage the spatiotemporal decomposition of the computational domain.
For instance, conservative PINNs (cPINNs)~\cite{jagtap2020conservative} utilize a non-overlapping Schwarz-based DD approach for conservation laws.
The extended PINNs (XPINNs) expand the cPINN methodology to generic PDEs and arbitrary space-time domains.
The idea behind XPINNs~\cite{jagtap2021extended,jagtap2020conservative, shukla2021parallel,hu2021extended} is to associate each subdomain with a sub-PINN, which can be trained independently of each other, i.e., in parallel.
After the training of the subdomains is concluded, the sub-PINNs are stitched together by enforcing continuity on the interfaces and by averaging solutions along neighboring subdomain interfaces.
The XPINN framework was later expanded in several ways.
For instance, the adaptive PINNs (APPINs)~\cite{hu2022augmented}, GatedPINNs~\cite{stiller2020large}, or extensions incorporating extreme learning machines~\cite{dwivedi2021distributed} have been proposed.

An alternative approach that utilizes the variational principle and overlapping Schwarz-based DD approach has also been explored, see~\cite{li2020deep}.
This method was further extended by introducing coarse space acceleration to improve convergence with respect to an increasing number of subdomains~\cite{mercier2021coarse}.
Building upon the variational principles and the deep Ritz method, the D3M~\cite{li2019d3m} also employs the overlapping Schwarz method.
\newtextreviewerone{A different approach was followed in~\cite{lee2021partition, trask2022hierarchical}, where the authors proposed partition of unity networks, which employ elements of domain decomposition to construct the localized basis functions to enhance the network's approximation properties.}
\oldtext{Another approach, called} \newtextreviewerone{Along similar lines, the} finite basis PINNs (FBPINN) \newtextreviewerone{approach, proposed in}~\cite{moseley2021finite, dolean2022finite}, represents the solution of a PDE as a sum of basis functions with compact support defined over small, overlapping subdomains.
\newtextreviewerone{Here, however the} \oldtext{These} basis functions are learned in parallel, \oldtext{using} \newtextreviewerone{by training a} neural network\oldtext{s} \newtextreviewerone{per subdomain}.
Finally, we mention PFNN-2 method~\cite{sheng2022pfnn}, which enhances the training of PINNs by utilizing overlapping DD together with the penalty-free network approach (PFNN)~\cite{sheng2021pfnn}.

In this work, we aim to enhance the training of the PINNs by developing novel nonlinear additive and multiplicative preconditioning strategies.
More specifically, we utilize the right preconditioning approach to accelerate the convergence of the widely used L-BFGS method.
The devised preconditioners are constructed by decomposing the network parameters in a layer-wise manner.
\newtext{Methodologically, our approach differs from the aforementioned DD approaches for PINNs, as we decompose the parameter space, while the decomposition of the computational domain is not considered.}
\oldtext{Methodologically, our approach differs from the aforementioned DD approaches for PINNs that utilize the decomposition of the computational domain.}
Consequently, the proposed preconditioners are agnostic to the physics-based priors and can potentially be applied to other machine-learning tasks.
Using the proposed DD approach, we demonstrate that both additive and multiplicative preconditioners can significantly improve the convergence rate of the standard L-BFGS optimizer.
Moreover, utilizing the additive preconditioner naturally enables model parallelism, which allows for additional acceleration of the training if a larger amount of computational resources is available.

This paper is organized as follows.
\Cref{sec:pinn} provides a brief overview of the PINNs.
In \cref{sec:preconditioner}, we propose nonlinear additive and multiplicative preconditioners using a layer-wise decomposition of the network.
\Cref{sec:implementation} discusses implementation details, while \cref{sec:experiments} introduces a set of benchmark problems used for numerical experiments.
\Cref{sec:results} demonstrates the numerical performance of the proposed preconditioning strategies.
Finally, \cref{sec:conclusion} summarizes the presented work.
 \section{Physics-informed neural networks}
\label{sec:pinn}
We consider the differential equations of a generic form and pose the problem as follows.
Find $u:\Omega \rightarrow \mathbb{R}$, such that
\begin{equation}
 \label{eq:pde_def}
 \begin{aligned}
  \pazocal{P}(u(\xv)) &= f(\xv), \quad   & \xv \in \Omega,               & &                         \\
  \pazocal{B}^k(u(\xv)) &= g^k(\xv), \quad & \xv\in \Gamma^k \subseteq \partial \Omega, & & \qquad  \text{for} \ k = 1,2,\ldots,n_\Gamma.
 \end{aligned}
\end{equation}
The symbol $\pazocal{P}$ represents a nonlinear differential operator, while $\{ \pazocal{B}^k\}_{k=1}^{n_\Gamma}$ denotes a set of boundary condition operators.
In the case of time-dependent problems, we treat time $t$ as an additional component of the vector $\xv \in \R^d$, where $d$ corresponds dimension of the problem.
Therefore, the computational domain $\Omega$ encompasses both spatial and temporal dimensions, and the initial conditions can be considered as a specific type of boundary conditions within the spatiotemporal domain.

\subsection{DNN approximation of solution}
Following~\cite{raissi2019physics}, our goal is to approximate a solution of~\eqref{eq:pde_def} using PINNs, i.e., ${u(\xv) \approx u_\texttt{NN}(\boldsymbol{\theta}, \xv)}$, where the nonlinear mapping ${u_\texttt{NN}: \R^{n} \times \R^d \rightarrow \R}$ is parametrized by the trainable parameters~${\boldsymbol{\theta} \in \R^n}$.
Typically, the exact form of~$u_\texttt{NN}$ is induced by the network architecture.
Here, we consider a network that is constructed as a composition of an input layer, $L-1$ hidden layers, and an output layer, thus as
\begin{equation}
 \begin{aligned}
  \text{input layer:}  & \qquad \yv_0 & = & \ \boldsymbol{W}_0 \xv,                                \\
  \text{hidden layers:} & \qquad \yv_l & = & \ \pazocal{N}(\yv_{l-1}), \qquad \qquad \text{for} \ \ l = 1, \ldots, L-1, \\
  \text{output layer:} & \qquad \yv_L & = & \ \boldsymbol{W}_{L} \yv_{L-1} + \bv_L.
 \end{aligned}
\end{equation}
The input layer maps input features~$\xv \in \Omega \subset \R^d$ into the dimension of the hidden layers using weights~${\boldsymbol{W}_0 \in \R^{d \times nh}}$.
The $l^{th}$ hidden layer nonlinearly transforms the output of the $(l-1)^{th}$ layer $\yv_{l-1}$ into $\yv_l$ by means of nonlinear operator~${\pazocal{N}: \R^{nh} \rightarrow \R^{nh}}$.
The operator~$\pazocal{N}$ can take on multiple forms.
We utilize a single-layer perceptron with a skip connection, given as
\begin{equation}
 \pazocal{N}(\yv_{l-1}) := \yv_{l-1} + \sigma_l(\boldsymbol{W}_l \yv_{l-1} + \bv_l).
\end{equation}
The parameters associated with layer~$l$ can be given as $\boldsymbol{\theta}_l := (\text{flat}(\boldsymbol{W}_l), \text{flat}(\bv_l) )$ collectively, where~$\boldsymbol{W}_l \in \R^{nh \times nh}$ and~$\bv_l \in \R^{nh}$ denote weights and biases, respectively.
The function~$\text{flat}(\cdot)$ represents a flattening operator and~$\sigma_l: \R^{nh} \rightarrow \R^{nh}$ denotes an activation function.
Finally, the output of the last hidden-layer~$\yv_{L-1}$ is linearly transformed using $\boldsymbol{W}_{L} $ and $\bv_L$, giving rise to the network's output~$\yv_{L}$.

\subsection{Training data and loss functional}
\label{sec:training_data_loss}
In order to train the network such that it can approximate a solution of~\eqref{eq:pde_def}, we construct a dataset~$\pazocal{D}_{\texttt{int}} =\{ \xv_j \}_{j=1}^{n_{\text{int}}}$ of $n_{\text{int}}$ collocations points, which lie in the interior of the computational domain~$\Omega$.
These points are used to enforce the physics captured by the PDE.
In addition, we consider datasets~$\{\pazocal{D}^k_{\texttt{bc}}\}_{k=1}^{n_\Gamma}$, where each ${\pazocal{D}^k_{\texttt{bc}} =\{ (\xv^k_j, \gv^k_j) \}_{j=1}^{n^k_{\text{bc}}}}$ contains~$n^k_{\text{bc}}$ points, which we use to enforce the boundary conditions on $\Gamma^k$.

Using $\pazocal{D}_{\texttt{int}}$ and $\{ \pazocal{D}^k_{\texttt{bc}}\}_{k=1}^{n_\Gamma}$, the optimal parameters of the network are found by solving the following minimization problem:
\begin{equation}
  \begin{aligned}
    {\boldsymbol{\theta}^\ast := \argmin_{\boldsymbol{\theta} \in \R^n} }
    \underbrace{ \tilde{\pazocal{L}}_{\text{int}} (\boldsymbol \theta)}_{\text{interior loss}} + \underbrace{ \tilde{\pazocal{L}}_{\text{bc}}(\boldsymbol{\theta})}_{\text{boundary loss}}.
  \end{aligned}
  \label{eq:min_problem_no_bc}
\end{equation}
Here, the interior and the boundary losses are given as
\begin{equation}
  \begin{aligned}
     & \tilde{\pazocal{L}}_{\text{int}}(\boldsymbol{\theta})  := \frac{w_{\text{int}}}{| \pazocal{D}_{\texttt{int}}|} \sum_{\xv_j \in \pazocal{D}_{\texttt{int}}} \| \pazocal{P} (\tilde{u}_\texttt{NN}(\boldsymbol{\theta}, \xv_j)\big) - f(\xv_j) \|^2,                                       \\
     & \tilde{\pazocal{L}}_{\text{bc}}(\boldsymbol{\theta}) := \sum_{k=1}^{n_\Gamma} \Bigg( \frac{w_{\text{bc}} }{| \pazocal{D}^k_{\texttt{bc}}|} \sum_{(\xv^k_j, \gv^k_j) \in \pazocal{D}^k_{\texttt{bc}}} \| (\tilde{u}_\texttt{NN}(\boldsymbol{\theta}, \xv^k_j) - g^k(\xv^k_j) \|^2 \Bigg).
  \end{aligned}
  \label{eq:loss}
\end{equation}
Thus, the interior loss~$\tilde{\pazocal{L}}_{\text{int}}$ is expressed as the mean PDE residual obtained by Monte Carlo integration, i.e., by the averaging over the residuals evaluated at collocation points sampled in the interior of the domain.
Note that an evaluation of the PDE residual typically requires knowledge of the partial derivatives of the network output~$\tilde{u}_\texttt{NN}(\boldsymbol{\theta}, \xv_j)$ with respect to the input~$\xv_j$.
These derivatives can be efficiently obtained using advanced automatic differentiation techniques.
The second term, boundary loss~$\tilde{\pazocal{L}}_{\text{bc}}$, ensures that the boundary conditions are satisfied at a set of points sampled along boundaries $\{\Gamma^k\}_{k=1}^{n_{\Gamma}}$.

The weights~$w_{\text{int}}$ and $w_{\text{bc}}$ in~\eqref{eq:loss} are used to balance the interior and boundary loss terms.
Their optimal value can be determined either during the hyper-parameter tuning or directly during the training by employing adaptive strategies, such as one proposed in~\cite{mcclenny2020self}.
\oldtext{In this work,} \newtextreviewerone{Following~\cite{lagaris1998artificial},} we overcome the difficulty of finding suitable values of~$w_{\text{int}}$ and $w_{\text{bc}}$ by imposing the boundary conditions in a penalty-free manner, thus by eliminating the boundary loss term. \oldtext{~\cite{lu2021physics}}
To this aim, the output of the network~$\tilde{u}_{\texttt{NN}}(\boldsymbol{\theta}, \xv)$ is modified using the length factor function~${\ell^k: \Gamma^k \subset \R^d \rightarrow \R}$ as follows
\begin{equation}
  \label{eq:mod_bc}
  {u}_{\texttt{NN}}(\boldsymbol{\theta}, \xv)  :=  \sum_{k=1}^{n_\Gamma} g^k (\xv)  +  \frac{\tilde{\ell}(\xv)}{\max\limits_{\xv \in \Omega} \tilde{\ell}(\xv)} \tilde{u}_\texttt{NN}(\boldsymbol{\theta}, \xv), \ \ \text{where } \tilde{\ell}(\xv) = \prod_{k=1}^{n_\Gamma} (1 - (1-\ell^k(\xv))).
\end{equation}
The function~$\ell^k$ is constructed such that
\begin{equation}
  \begin{aligned}
    \ell^k(\xv) & =  0, \quad      &  & \xv \in \Gamma^k                                                                    \\
\ell^k(\xv) & =  1, \quad      &  & \xv \in \Gamma^{j},\quad j\in\{1,2,\ldots,n_\Gamma\},  \quad \text{for} \  k \neq j \\
    \ell^k(\xv) & \in (0,1)  \quad &  & \text{otherwise.}                                                                   \\
  \end{aligned}
\end{equation}
The details regarding the construction of~$\ell^k$ for complex geometries can be found in~\cite{sukumar2022exact}.

Using~\eqref{eq:mod_bc}, the minimization problem~\eqref{eq:min_problem_no_bc} can be now reformulated as
\begin{equation}
  \boldsymbol{\theta}^\ast := \argmin_{\boldsymbol{\theta} \in \R^n} \ \pazocal{L}(\boldsymbol{\theta}) :=\frac{1}{| \pazocal{D}_{\texttt{int}}|}  \sum_{\xv_j \in \pazocal{D}_{\texttt{int}}} \| \pazocal{P} ({u}_\texttt{NN}(\boldsymbol{\theta}, \xv_j) \big) - f(\xv_j) \|^2.
  \label{eq:min_problem}
\end{equation}
After training, the optimal parameters~$\boldsymbol{\theta}^\ast$ can be used to infer an approximate solution~$u_{\texttt{NN}}$ of~\eqref{eq:pde_def}.
The quality of this solution can be assessed by means of the error $\pazocal{E}(u_\texttt{NN}, u^*):= \| u_{\texttt{NN}} - u^* \|$, where $u^*$ denotes an exact solution.
In general, the error~$\pazocal{E}$ consists of three components: the discretization error, network approximation error, and optimization error~\cite{mishra2022estimates}.
More specifically, the discretization error is determined by the number and location of the collocation points.
The approximation error of the network is associated with the representation capacity of the DNN, i.e., it directly depends on the network architecture (number of layers, network width).
The optimization error quantifies the quality of the solution provided by the optimizer.
Our goal in this work is to reduce the optimization error by developing enhanced training algorithms.

\section{Nonlinearly preconditioned training using layer-wise decomposition of network}
\label{sec:preconditioner}
The minimization problem~\eqref{eq:min_problem} is traditionally solved using Adam and/or L-BFGS optimizers.
As it has been pointed out in~\cite{wang2021understanding}, the convergence of these optimizers deteriorates with an increasing problem stiffness, which is caused by the multiscale and multirate nature of the considered problems.
We aim to overcome this difficulty by proposing a novel nonlinear preconditioning strategy.

Let us recall that a minimizer of~\eqref{eq:min_problem} has to satisfy the following first-order criticality condition
\begin{equation}
	\nabla \pazocal{L}(\boldsymbol{\theta}) = 0.
	\label{eq:res_eq_zero}
\end{equation}
As solving~\eqref{eq:res_eq_zero} is computationally demanding, we instead construct and solve an alternative, a nonlinearly preconditioned system of equations, given as
\begin{equation}
	\pazocal{F}(\boldsymbol{\theta})=0.
	\label{eq:precond_system}
\end{equation}
The preconditioned system of equations~\eqref{eq:precond_system} shall be constructed such that it is easier to solve than~\eqref{eq:res_eq_zero}, but such that it has the same solution as an original nonlinear system.
Thus, the preconditioned problem should have more balanced nonlinearities, and the numerical computations involved in solving \eqref{eq:precond_system} should be computationally cheaper.

\subsection{Network decomposition}
\label{sec:net_decomposition}

\begin{figure}
	\centering
	\def\layersep{1.8 cm}
	\def\center_y{-1.5cm}
\includegraphics{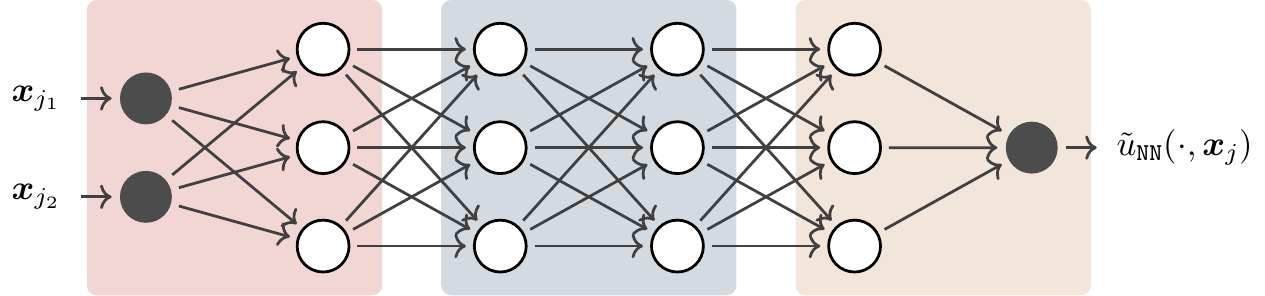}

\caption{An example of the  layer-wise decomposition of the network.}
	\label{fig:decomp_net}
\end{figure}
 In this work, we propose to construct the nonlinearly preconditioned system of equations~\eqref{eq:precond_system} by decomposing the parameter space of the network.
In particular, subnetworks are constructed by layer-wise decomposition of the network into~$N_{sd}$ disjoint groups, often called subdomains.
Thus, each layer is assigned a unique subnetwork identifier, and the parameters of the network are also split into~$N_{sd}$ disjoint groups, i.e.,~$\boldsymbol{\theta} = [ \boldsymbol{\theta}_1, \ldots, \boldsymbol{\theta}_s, \ldots, \boldsymbol{\theta}_{N_{sd}} ]^\top $, where~$\boldsymbol{\theta}_s \in \R^{n_{s}}$ denotes parameters associated with a $s^{th}$ subnetwork.
Please refer to \Cref{fig:decomp_net} for the illustration of layer-wise decomposition.
\oldtext{In passing, we note that the different types of network decomposition can be potentially utilized, e.g., intra-layer decomposition with or without overlap.}
\newtextreviewerone{Here, we point out that the different types of network decompositions can be utilized to construct the nonlinear preconditioner, e.g., intra-layer decomposition with or without overlap.}

Transfer of the data between the subnetwork and the global network is performed using two transfer operators.
Restriction operator $\boldsymbol{R}_s: \R^n \rightarrow \R^{n_{s}}$ extracts the parameters of network associated with the $s^{th}$ subnetwork, i.e.,
\begin{equation}
	\boldsymbol{\theta}_s = \boldsymbol{R}_s \boldsymbol{\theta}, \qquad \text{for} \ s = 1, 2, \ldots, N_{sd}.
\end{equation}
In contrast, an extension operator $\boldsymbol{E}_s: \R^{n_{s}} \rightarrow \R^{n}$ is used to extend the quantities from the  subnetworks to the global network, i.e.,
\begin{equation}
	\boldsymbol{\theta} = \sum_{s=1}^{N_{sd}} \boldsymbol{E}_s \boldsymbol{\theta}_s.
\end{equation}

Using the proposed network decomposition, we now define local minimization problems associated with a given subnetwork~$s$ as
\begin{equation}
	\boldsymbol{\theta}_s^* = \argmin_{\boldsymbol{\theta}_s} \pazocal{L} (\boldsymbol{\theta}_1, \ldots, \boldsymbol{\theta}_s, \ldots, \boldsymbol{\theta}_{N_{sd}}).
	\label{eq:min_local}
\end{equation}
Solving~\eqref{eq:min_local} corresponds to minimizing the loss functional~$\pazocal{L}$ with respect to the parameters~$\boldsymbol{\theta}_s$, while all other parameters of the network are kept fixed.
This minimization process can be carried out by any optimizer of a choice, such as Adam~\cite{kingma2014adam}, L-BFGS~\cite{liu1989limited}, or Newton~\cite{deuflhard2011newton}.
For the purpose of this work, we employ the L-BFGS method, see also \cref{sec:implementation}.
Note, this process loosely resembles the transfer-learning approaches~\cite{weiss2016survey}, where only a particular subset of the network parameters is being (re-)trained.

\subsection{Right-preconditioned training optimizers}
In this work, we construct a nonlinearly preconditioned system of equations as follows
\begin{equation}
	\pazocal{F}(\boldsymbol{\theta}) := \nabla \pazocal{L}(G(\boldsymbol{\theta})),
	\label{eq:composite_precond}
\end{equation}
where the operator~$G:\R^{n} \rightarrow \R^n$ denotes a right preconditioner.
Using the network-decomposition framework introduced in \cref{sec:net_decomposition}, we define operator~$G$ as an iterative method of the following form:
\begin{equation}
	\newtext{ \boldsymbol{\theta}^{(k+1/2)} =:} G(\boldsymbol{\theta}^{(k)}) = \boldsymbol{\theta}^{(k)} + \alpha^{(k)} \sum_{s=1}^{N_{sd}} \boldsymbol{E}_s (\boldsymbol{\theta}_s^* - \boldsymbol{R}_s \boldsymbol{\theta}^{(k)}).
\end{equation}
Here, a symbol~$\boldsymbol{\theta}_s^*$ represents a solution of the local minimization problem~\eqref{eq:min_local}, associated with $s^{th}$ subnetwork, while~$\alpha^{(k)} $ denotes a step-size.

\oldtext{The subnetwork solutions $\{\boldsymbol{\theta}_s^*\}_{s=1}^{ N_{sd}}$ can be obtained independently from each other, i.e., in parallel, which gives rise to an additive preconditioner~$G$.
Alternatively, one might obtain $\{\boldsymbol{\theta}_s^*\}_{s=1}^{ N_{sd}}$ in a multiplicative fashion.
In this case, the local minimization problems are solved sequentially and the already acquired local solutions $\{\boldsymbol{\theta}_z^*\}_{z=1}^{s-1}$ are taken into account while performing local training for the $s^{th}$ subnetwork.}
\newtext{The subnetwork solutions $\{\boldsymbol{\theta}_s^*\}_{s=1}^{ N_{sd}}$ can be obtained in an additive or a multiplicative manner.
In the additive case, the local minimization problems are solved in parallel, starting from~$\boldsymbol{R}_s \boldsymbol{\theta}^{(k)}$.
In the multiplicative case, the local minimization problems are solved sequentially, i.e., one after another.
As soon as the solution of the $s^{th}$ local problem is obtained, it is taken into account while performing local training for the $(s+1)^{th}$ subnetwork.
In this way, the multiplicative approach allows for faster information transfer between the different subnetworks, albeit at the price of obtaining a sequential algorithm.}

As we can see from~\eqref{eq:composite_precond}, $\pazocal{F}$ is a composite operator.
Thus, for a given iterate~$\boldsymbol{\theta}^{(k)}$, the preconditioner~$G$ is used to obtain an improved iterate, denoted by~${\boldsymbol{\theta}^{(k+\sfrac{1}{2})}= G(\boldsymbol{\theta}^{(k)})}$.
The improved iterate~$\boldsymbol{\theta}^{(k+\sfrac{1}{2})}$ is then used as an initial guess for the next global optimization step,\footnote{\newtext{The terminology ``global'' used in this work refers to quantities associated with the global network (the original network), and it differs from the concept of ``global optimization'', which is concerned with finding the global minimum.}}
giving rise to the following update rule:
\begin{equation}
	\boldsymbol{\theta}^{(k+1)} = \boldsymbol{\theta}^{(k+\sfrac{1}{2})} + \alpha^{(k+\sfrac{1}{2})} \pv^{(k+\sfrac{1}{2})},
\end{equation}
\oldtext{where~$\pv^{(k+\sfrac{1}{2})}$ is obtained by a global optimizer and~$\alpha^{(k+\sfrac{1}{2})}$ denotes appropriately selected global step-size.}
\newtext{where~$\alpha^{(k+\sfrac{1}{2})}$ is an appropriately selected global step-size, while~$\pv^{(k+\sfrac{1}{2})}$ is a search direction obtained by performing a single step of optimizer for the global network.
For example, if global optimizer is chosen to be a gradient descent method, then search direction would be obtained as $\pv^{(k+\sfrac{1}{2})}:= - \nabla \pazocal{L} (\boldsymbol{\theta}^{(k+\sfrac{1}{2})})$.}

\subsubsection{Right-preconditioned L-BFGS}
There are many possible choices of a global optimizer, e.g., Adam~\cite{kingma2014adam} or AdaGrad~\cite{lydia2019adagrad}.
Here, we utilize the BFGS method, as it is widely employed for the training of PINNs.
The BFGS method belongs to the family of quasi-Newton (QN) methods, which have gained popularity due to the fact that they do not require explicit knowledge of the Hessian.
Instead, a Hessian approximation~$\boldsymbol{B}^{(k+1)}$ is constructed such that the following secant equation is satisfied:
\begin{equation}
	\boldsymbol{B}^{(k+1)} \sv^{(k)} = \yv^{(k)},
	\label{eq:secant_eq}
\end{equation}
where $\sv^{(k)} = \boldsymbol{\theta}^{(k+1)} - \boldsymbol{\theta}^{(k+\sfrac{1}{2})}$ and $\yv^{(k)} = \nabla \pazocal{L}(\boldsymbol{\theta}^{(k+1)}) - \nabla\pazocal{L}(\boldsymbol{\theta}^{(k+\sfrac{1}{2})})$.
Note that in the context of right-preconditioned nonlinear systems, the secant pairs are obtained by taking into account iterate~$\boldsymbol{\theta}^{(k+\sfrac{1}{2})}$, obtained as an outcome of the nonlinear preconditioning step, see~\cite{kothari2022nonlinear} for details.
Since the secant equation~\eqref{eq:secant_eq} is indeterminate, additional constraints must be applied to uniquely determine the approximation~$\boldsymbol{B}^{(k+1)}$.
A particular choice of constraints leads to different variants of the QN methods.
For instance, the BFGS method is obtained by imposing the symmetry and positive definiteness of $\boldsymbol{B}^{(k+1)}$.

In order to reduce the memory footprint of the proposed training method, we employ a limited-memory variant of the BFGS method, termed L-BFGS.
Using compact matrix representation~\cite{liu1989limited}, the approximation~$\boldsymbol{B}^{(k+1)}$ is given as
\begin{equation}
\boldsymbol{B}^{(k+1)} =
	\boldsymbol{B}^{(0)} -
	\begin{bmatrix} \boldsymbol{B}^{(0)}  \boldsymbol{S}^{(k)}  \\ \boldsymbol{Y}^{(k)}\end{bmatrix}^\top
	\begin{bmatrix} (\boldsymbol{S}^{(k)})^\top \boldsymbol{B}^{(0)} \boldsymbol{S}^{(k)} & \Lm^{(k)} \\ (\Lm^{(k)})^\top & - \Dm^{(k)} \end{bmatrix}^{-1}
	\begin{bmatrix} (\boldsymbol{S}^{(k)})^\top \boldsymbol{B}^{(0)} \\ (\boldsymbol{Y}^{(k)})^\top\end{bmatrix},
	\label{eq:LBFGS}
\end{equation}
where the matrices $\boldsymbol{S}^{(k)} := [\sv^{(k-m)},\ldots, \sv^{(k-1)}]$, and $\boldsymbol{Y}^{(k)} := [\yv^{(k-m)}, \ldots, \yv^{(k-1)}]$ are constructed using $m$ most recent secant pairs~$\{(\sv^{(i)}, \yv^{(i)}) \}_{i=k-m}^{k-1}$.
The matrices $\Lm^{(k)}$ and $\Dm^{(k)}$ denote the strictly lower triangular, and the diagonal part of the matrix $(\boldsymbol{S}^{(k)})^\top \boldsymbol{Y}^{(k)}$, while $\boldsymbol{B}^{(0)}$ is an initial approximation of the Hessian.
Usually $\boldsymbol{B}^{(0)} = \gamma \Im$, where $\gamma = \frac{\langle \yv^{(k)}, \yv^{(k)} \rangle}{\langle \yv^{(k)}, \sv^{(k)} \rangle} $,  see~\cite{rafati2018improving} for alternatives.

Using~\eqref{eq:LBFGS}, we are now able to state the global L-BFGS step as follows
\begin{equation}
	\boldsymbol{\theta}^{(k+1)} = \boldsymbol{\theta}^{(k+\sfrac{1}{2})} + \alpha^{(k+\sfrac{1}{2})} \underbrace{\big(-( \boldsymbol{B}^{(k+1)} )^{-1} \nabla \pazocal{L}(\boldsymbol{\theta}^{(k+\sfrac{1}{2})})\big)}_{=:\pv^{(k+\sfrac{1}{2})}}.
	\label{eq:global_step_no_momentum}
\end{equation}
Using the compact matrix representation, the inverse Hessian~$( \boldsymbol{B}^{(k+1)} )^{-1}$ can be obtained explicitly by applying the Sherman--Morrison--Woodbury formula.
Alternatively, the two-loop recursion formula~\cite{gilbert1993automatic} can be utilized.

The convergence of the right-preconditioned L-BFGS method can be further enhanced by incorporating momentum, as often done for training DNNs.
Following~\cite{mahboubi2019momentum, rafati2018improving}, we evaluate the momentum~$\vv^{(k+\sfrac{1}{2})} \in \R^{n}$ in recursive manner as
\begin{equation}
	\vv^{(k+\sfrac{1}{2})} = (1-\mu) \vv^{(k-\sfrac{1}{2})} + \mu \pv^{(k+\sfrac{1}{2})},
\end{equation}
where~$\mu \in [0,1]$.
Utilizing the momentum, the global step~\eqref{eq:global_step_no_momentum} can be now reformulated as
\begin{equation}
	\boldsymbol{\theta}^{(k+1)} = \boldsymbol{\theta}^{(k+\sfrac{1}{2})} + \alpha^{(k+\sfrac{1}{2})} \vv^{(k+\sfrac{1}{2})}.
\end{equation}

\Cref{alg:rplbfgs} summarizes the Schwarz preconditioned quasi-Newton (SPQN) methods, namely additive and multiplicative variants termed ASPQN and MSPQN, respectively.

\begin{algorithm}[t]
	\caption{Schwarz Preconditioned Quasi-Newton (SPQN)}
	\label{alg:rplbfgs}
	\begin{algorithmic}[1]
		\REQUIRE $\pazocal{L}: \R^n \rightarrow \R, \  \boldsymbol{\theta}^{(0)} \in \R^n, \  \epsilon \in \R, \  \mu \in [0,1],  \ k_{\text{max}} \in \N^+, \vv^{-\sfrac{1}{2}} =
			\boldsymbol{0}$
		\STATE  $ \vv^{(0)} = 0 $
		\FOR{$k = 0, \ldots, k_{\text{max}}$ }
		\FOR{$s = 1, \ldots, N_{sd}$}
		\IF{additive}
		\STATE $\boldsymbol{\theta}^{*}_s  = \argmin_{\boldsymbol{\theta}_s} \pazocal{L}(\boldsymbol{\theta}_1^{(k)}, \ldots,  \boldsymbol{\theta}_s, \ldots, \boldsymbol{\theta}_{N_{sd}}^{(k)})$
		\COMMENTmine{parallel}
		\ELSE \STATE $\boldsymbol{\theta}^{*}_s  = \argmin_{\boldsymbol{\theta}_s} \pazocal{L}(\boldsymbol{\theta}^*_1, \ldots,  \boldsymbol{\theta}^*_{s-1}, \boldsymbol{\theta}_s, \boldsymbol{\theta}^{(k)}_{s+1}, \ldots, \boldsymbol{\theta}_{N_{sd}}^{(k)})$
		\COMMENTmine{sequential}
		\ENDIF
		\ENDFOR
		\STATE
		\STATE $ \boldsymbol{\theta}^{(k+\sfrac{1}{2})} =  \boldsymbol{\theta}^{(k)} + \alpha^{(k)} \sum_{s=1}^{N_{sd}} \boldsymbol{E}_s (\boldsymbol{R}_s \boldsymbol{\theta}^{(k)} - \boldsymbol{\theta}_s^*) $   \COMMENTmine{Synchronization step}
		\STATE
		\STATE $\pv^{(k+\sfrac{1}{2})} = - (\boldsymbol{B}^{(k+1)})^{-1} \nabla \pazocal{L}(\boldsymbol{\theta}^{(k+\sfrac{1}{2})}) $   \COMMENTmine{Preconditioned quasi-Newton step}
		\STATE $\vv^{(k+\sfrac{1}{2})} = (1-\mu) \vv^{(k-\sfrac{1}{2})} +  \mu \pv^{(k+\sfrac{1}{2})} $   \COMMENTmine{Momentum evaluation}
		\STATE $ \boldsymbol{\theta}^{(k+1)} =\boldsymbol{\theta}^{(k+\sfrac{1}{2})} + \alpha^{(k+\sfrac{1}{2})} \vv^{(k+\sfrac{1}{2})}$  \COMMENTmine{Global parameter update}
		\STATE
		\STATE Update~$\boldsymbol{S}^{(k+1)}$ with $\sv^{(k)}= \boldsymbol{\theta}^{(k+1)} - \boldsymbol{\theta}^{(k+\sfrac{1}{2})}  $ \COMMENTmine{Update of secant pairs}
		\STATE Update~$\boldsymbol{Y}^{(k+1)}$ with $\yv^{(k)}= \nabla \pazocal{L}(\boldsymbol{\theta}^{(k+1)}) -  \nabla \pazocal{L}(\boldsymbol{\theta}^{(k+\sfrac{1}{2})})  $
		\ENDFOR
		\RETURN $\boldsymbol{\theta}^{(k+1)}$
	\end{algorithmic}
\end{algorithm}
 \section{Numerical experiments}
\label{sec:experiments}
In this section, we describe numerical examples used to test the proposed SPQN methods.
In particular, we consider the following benchmark problems.
\begin{itemize}[leftmargin=0.5cm]
	\item{\textbf{Burgers' equation:}}
	      Burgers' equation has the following form:
	      \begin{equation}
		      \begin{aligned}
			      \frac{\partial u}{\partial t} + u \nabla u - \nu \nabla^2 u & = 0,  \quad \quad \quad &  & \forall \ (t,x) \in (0,1] \times (-1,1), \\
			      u                                                           & = - \sin(\pi x), \      &  & \forall \ (t,x) \in \{0\} \times [-1,1], \\
			      u                                                           & = 0, \                  &  & \forall \ (t,x) \in (0,1] \times \{1\},  \\
			      u                                                           & = 0, \                  &  & \forall \ (t,x) \in (0,1] \times \{-1\},
		      \end{aligned}
	      \end{equation}
	      where~$u=u(t,x)$ denotes the flow velocity and~$\nu$ stands for the kinematic viscosity.
	      Here, we choose~$\nu=0.01/\pi$.

	\item{\textbf{Diffusion-advection:}}
	      The steady-state diffusion-advection problem is given in two-spatial dimensions as
	      \begin{equation}
		      \begin{aligned}
			      -\nabla \cdot \mu  \nabla u  + \boldsymbol{b}\cdot \nabla u & = f, &  & \forall \ (x_1,x_2) \in (0, 1) \times (0, 1), \\
			      u                                                           & = 0, &  & \text{on } \partial \Omega,
		      \end{aligned}
	      \end{equation}
	      where~$\boldsymbol{b}=(1,1)^\top$.
	      The data is considered to be constant, i.e.,  $f = 1$.
	      Moreover, the symbol~$\mu$ denotes viscosity.
	      In this work, we employ $\mu=10^{-2}$.
	\item{\textbf{Klein-Gordon:}}
	      The nonlinear second-order hyperbolic equation, called Klein-Gordon, is defined as
	      \begin{equation}
		      \begin{aligned}
			      \frac{\partial^2 u}{\partial t^2}  + \alpha \nabla^2 u + \beta u + \gamma u^2 & = f(t,x), \quad \quad \quad &  & \forall \ (t,x) \in (0,12] \times (-1,1), \\
			      u                                                                             & = x, \                      &  & \forall \ (t,x) \in \{0\} \times [-1,1],  \\
			      \frac{\partial u}{\partial t}                                                 & = 0, \                      &  & \forall \ (t,x) \in \{0\} \times [-1,1],  \\
			      u                                                                             & = - \cos(t), \              &  & \forall \ (t,x) \in (0,12] \times \{-1\}, \\
			      u                                                                             & = \cos(t), \                &  & \forall \ (t,x) \in (0,12] \times \{1\},
		      \end{aligned}
	      \end{equation}
	      where~$u=u(t,x)$ denotes the wave displacement at time~$t$ and location~$x$.
	      The symbols~$\alpha, \beta, \gamma$ denote the predetermined constants.
	      For the purpose of this work, we employ~${\alpha=-1, \beta=0, \gamma=1}$ and ${f(t,x):= -x \cos(t) + x^2 \cos^2(t)}$.
	      Given this particular choice of the parameters, the analytical solution is given as~$u(t,x)=x \cos(t)$, see~\cite{guo2015numerical}.

	\item{\textbf{Allen-Cahn:}}
	      We consider the Allen-Cahn equation of the following form:
	      \begin{equation}
		      \begin{aligned}
			      \frac{\partial u}{\partial t}  - D \nabla^2 u - 5 (u - u^3) & = 0,  \quad \quad \quad &  & \forall \ (t,x) \in (0,1] \times (-1,1), \\
			      u                                                           & = x^2 \cos(\pi x), \    &  & \forall \ (t,x) \in \{0\} \times [-1,1], \\
			      u                                                           & = - 1, \                &  & \forall \ (t,x) \in (0,1] \times \{-1\}, \\
			      u                                                           & = -1, \                 &  & \forall \ (t,x) \in (0,1] \times \{1\},
		      \end{aligned}
	      \end{equation}
	      where the diffusion coefficient~$D$ is chosen as $D=0.001$.
\end{itemize}

\section{Implementation, computational cost, and memory requirements}
\label{sec:implementation}
All presented numerical experiments were performed at the Swiss National Supercomputing Centre (CSCS) using XC50 compute nodes of the Piz Daint supercomputer.
Each XC50 compute node comprises an Intel Xeon E5-2690 V3 processor and an NVIDIA Tesla P100 GPU.
The memory of a node is $64$~GB, while the memory of the GPU is $16$~GB.

\subsection{\newtext{Implementation of PINNs}}
\begin{figure}
	\centering
	\def\layersep{1.8 cm}
	\def\center_y{-1.5cm}
\includegraphics{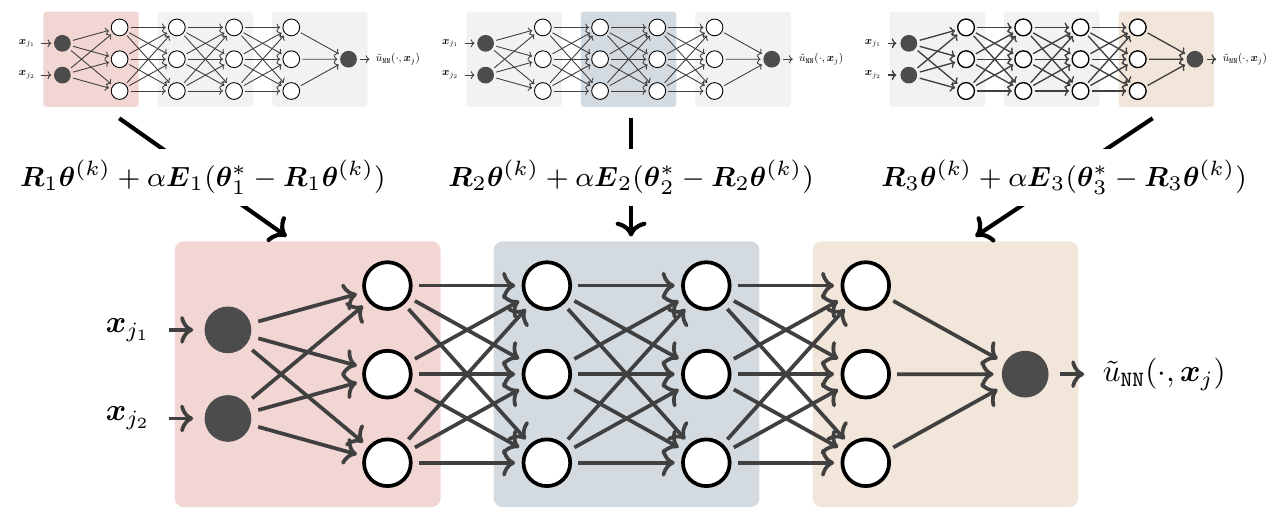}

\caption{A sketch of local-to-global updates utilized by the ASPQN.}
	\label{fig:parallelisation}
\end{figure}

\oldtext{We implement the proposed ASPQN and MSPQN optimizers as part of {\sf DistTraiNN} library~\cite{disttraingit}, which utilizes {\sf PyTorch}~\cite{paszke2019pytorch} for the implementation of PINNs.
As mentioned earlier, the MSPQN algorithm is inherently sequential, and therefore, its implementation targets computing architectures with a single GPU device.}
\newtext{We implement PINNs utilizing the {\sf PyTorch} library~\cite{paszke2019pytorch}.}
\newtext{The training set is constructed by using a uniform non-adaptive sampling strategy based on quasi-random Hammersley low-discrepancy sequences~\cite{mishra2022estimates}.}
The number of collocation points is chosen to be $10,000$ for all numerical examples.
\newtext{As discussed in \cref{sec:training_data_loss}, the error $\pazocal{E}(u_\texttt{NN}, u^*)$ consists of three components: the discretization error, network approximation error, and optimization error~\cite{mishra2022estimates}.
By choosing the number of collocation points deliberately high, we aim to reduce the discretization error as much as possible, which in turn allows us to study the impact of the optimization error in an isolated manner.}

\begin{table}[t]
  \centering
  \caption{Summary of network architecture and optimizer's parameters used for numerical experiments reported in \cref{sec:comparison}.}
  \label{tab:net_architectures}
\begin{tabular}{|l|l|l|l|l|l|l|l|l|l|l|l|l|l|l|l|l|l|}
    \hline
    \multirow{3}{*}{{Benchmark problem} } & \multicolumn{7}{c|}{{Parameter type}}                                                                                                                                    \\  \cline{2-8}
                                          & \multirow{1}{*}{{Depth }}            & \multirow{1}{*}{{Width}} & Adam                 & \multicolumn{2}{c|}{{ASPQN}} & \multicolumn{2}{c|}{{MSPQN}}                    \\ \cline{5-8}
                                          & $L$                                   & $nh$                     & $lr$                   & $N_{sd}$                     & $k_s$                        & $N_{sd}$ & $k_s$ \\ \hline \hline
    Burgers'                              & 8                                     & 20                       & $5 \times 10^{-4}$   & 8                            & 50                           & 8        &  50     \\  \hline
    Diffusion-advection                   & 10                                    & 50                       & $1 \times 10^{-4}$   & 10                           & 10                           & 10       & 10    \\  \hline
    Klein-Gordon                          & 6                                     & 50                       & $1 \times 10^{-3}$   & 6                            & 50                           & 6        &    50   \\  \hline
    Allen-Cahn                            & 6                                     & 64                       & $2.5 \times 10^{-4}$ & 6                            & 50                           & 6        & 50    \\  \hline
  \end{tabular}

\end{table}

\newtext{Similarly, we have selected network architectures such that their representational capacity is sufficiently large.}
\newtext{In particular, the network's width and depth, reported in \cref{tab:net_architectures}, were found by the hyperparameter search, performed using standard L-BFGS optimizer.
The width was selected from the set $\{20, 32, 50, 64, 100, 128 \}$ and depth from the set $\{ 4, 6, 8, 10, 12 \}$, such that they collectively gave rise to the lowest $\pazocal{E}(u_\texttt{NN}, u^*)$.}
As an activation function, we employ an adaptive variant of \textit{tanh}~\cite{jagtap2020adaptive}.
All networks are initialized using the Xavier initialization strategy~\cite{glorot2010understanding}.
\oldtext{To construct the training set, we use a uniform non-adaptive sampling strategy based on quasi-random Hammersley low-discrepancy sequences~\cite{mishra2022estimates}.}

\subsection{\newtext{Implementation of optimizers}}
The proposed ASPQN and MSPQN optimizers are implemented as part of the library {\sf DistTraiNN}\footnote{The code will be made publicly available upon acceptance of the manuscript.}~\cite{disttraingit}.
As mentioned earlier, the MSPQN algorithm is inherently sequential, and therefore, its implementation targets computing architectures with a single GPU device.
In contrast, the ASPQN algorithm is inherently parallel; its implementation is designed for distributed computing environments with multiple GPU devices.
Here, we utilize the {\sf torch.distributed} package with the NCCL backend and implement the ASPQN method such that each GPU is associated with a single subnetwork, i.e., the number of GPUs determines the number of subnetworks.
As an evaluation of loss associated with a subnetwork requires input to pass through an entire network, we replicate the global network to all GPUs.
Each subnetwork then undergoes local optimization by training the subnetwork's parameters, i.e., by minimizing~\eqref{eq:min_local}.
This is in contrast to other model parallel approaches, where each GPU gets access to only a specific part of the network~\cite{dean2012large, ben2019demystifying}.
As a consequence, the standard model parallel approaches have a smaller memory footprint compared to our approach.
After the subnetwork's training is concluded, the synchronization of local updates is performed (line 11, \Cref{alg:rplbfgs}).
Subsequently, a global L-BFGS step is executed concurrently on all nodes.
\Cref{fig:parallelisation} provides a sketch of the ASPQN parallelization strategy.

\begin{remark}
	The implementation of the SPQN methods can be seamlessly integrated with data parallel approaches.
\end{remark}

\subsubsection{\newtext{Configuration of SPQN optimizers}}
We solve the local minimization problem~\eqref{eq:min_local} using the L-BFGS method.
This minimization is carried out only approximately, using a fixed number of iterations, denoted by $k_s$.
Both, global and subnetwork, L-BFGS Hessian approximations take advantage of three secant pairs, i.e., $m=3$.
The Hessian approximation for each subnetwork is restarted every time the local training commences.
Moreover, we employ cubic backtracking line-search method with strong Wolfe conditions~\cite[Algorithm A6.3.1, pages 325-327]{dennis1996numerical} to obtain $\alpha^{(k+\sfrac{1}{2})}$ and $\alpha^{(k)}$.
Alternatively, one could consider using fixed step sizes, which may result in a more efficient algorithm that would require fewer evaluations of the loss function.
However, determining suitable values of $\alpha^{(k+\sfrac{1}{2})}$ and $\alpha^{(k)}$ would necessitate a thorough hyper-parameter search, which would eventually lead to a significant rise in overall computational cost.
\newtext{Unless specified differently, the SPQN methods are configured as summarized in \Cref{tab:net_architectures}.
As we will see in \cref{sec:experiments_decomp_study}, these particular configurations of SPQN methods give rise to the best-performing variants.}

\subsubsection{\newtext{Configuration of state-of-the-art optimizers}}
\newtext{In order to assess the performance of the proposed SPQN methods, we compare their performance with the standard Adam and L-BFGS optimizers.
The learning rate for the Adam optimizer is chosen by carrying out a thorough hyper-parameter search.
In particular, we choose a learning rate (reported in \Cref{tab:net_architectures}), which allows us to obtain the most accurate solution of the PDE, thus a solution with the lowest~$\pazocal{E}_{\text{rel}}$.
For the L-BFGS optimizer, we employ the same memory size ($m=3$), line-search strategy, and momentum setup as for the L-BFGS methods employed within the SPQN framework.
To ensure a fair comparison, Adam and L-BFGS optimizers are executed in deterministic settings, i.e., evaluation of loss/derivatives takes into account all samples in the dataset.}

\subsection{Computational cost and memory requirements}
In this section, we discuss the computational cost and memory requirements of the proposed SPQN methods, as well as standard Adam and L-BFGS optimizers.
\newtext{The summary is provided in } \cref{tab:computational_cost}, \newtext{where we compare the number of the loss function evaluations $(\# \pazocal{L}_e)$, the number of gradient evaluations  $(\# g_e)$, the update cost (UC) and memory cost (MC) per iteration and a device. }

\oldtext{In particular, the Adam optimizer requires} \newtext{In the case of the Adam optimizer, we require} a single evaluation of the loss function and the gradient per iteration.
The \oldtext{update cost} \newtext{UC} involves not only updating the parameters but also computing the biased and corrected first-order and second-order moment estimates, thus approximately $5n$ flops.
The memory requirement for the Adam optimizer is related to storing the iterate, gradient, and first-order and second-order moment estimates, leading to the total memory of size $4n$.

The L-BFGS method requires an evaluation of one gradient per iteration.
However, the loss function might be evaluated multiple times, depending on the number of required line-search iterations, denoted by~$ls_{\text{its}}$.
The update cost depends on the number of stored secant pairs ($m$) that are used for approximating the Hessian.
Using compact matrix representation, it can be approximately estimated as $4mn$ flops,
see~\cite[Sections 3 and 4]{byrd1994representations} for details.
Moreover, we have to take into account the cost associated with an evaluation of the momentum term, leading to an overall cost of approximately~$n+4mn$ flops.
In the context of memory requirements, the L-BFGS method necessitates the storage of the momentum term~$\vv^{(k)}$ as well as matrices~$\boldsymbol{S}^{(k)}$ and $\boldsymbol{Y}^{(k)}$, which accounts for the memory of size $n+2mn$.

In the case of the SPQN methods, the computational cost consists of global and local parts.
The global UC and memory requirements consist of that of the standard {L-BFGS} method.
Moreover, one has to take into account the cost related to the evaluation of the preconditioned iterate~$\boldsymbol{\theta}^{(k+\sfrac{1}{2})}$ (line 11, \Cref{alg:rplbfgs}).
The local UC is associated with the cost required by the local optimizer, denoted by $\text{UC}_s$.
Thus, for the $s^{th}$ subnetwork, \oldtext{ the update cost is scaled proportionally to the number of local parameters~$n_s$ and the number of local iterations~$k_s$.}
\newtext{we scale the $\text{UC}_s$ proportionally to the number of local parameters~$n_s$ and the number of local iterations~$k_s$.}
Assuming that all subnetworks consist of the same number of local trainable parameters ($n_s$), the local update cost of the ASPQN method is therefore~$(\sfrac{n_s}{n} )k_s  \text{UC}_s$.
For the MSPQN method, \oldtext{the local update cost has to be summed up over all subnetworks,} \newtext{we sum up the local update cost over all subnetworks} due to the sequential nature of the algorithm, thus $\sum_{s=1}^{N_{sd}} \big((\sfrac{n_s}{n}) k_s \text{UC}_s \big) = k_s \text{UC}_s$.
Regarding the \oldtext{memory cost} \newtext{MC} of both SPQN methods, we have to take into account the storage of global and local Hessian approximations, the momentum, and the preconditioned iterate~$\boldsymbol{\theta}^{(k+\sfrac{1}{2})}$.

In contrast to the update cost, \oldtext{the estimate for} \newtext{we do not scale} the number of the loss and gradient evaluations required by SPQN methods \oldtext{is not scaled} proportionally to $n_s$.
This is due to the fact that evaluation of the loss requires a forward pass through the entire network.
In contrast, the evaluation of the local gradient, denoted as $\nabla_{\boldsymbol{\theta}_s} \pazocal{L}$, does not necessitate traversing the entire network as several code optimization strategies can be utilized to avoid unnecessary computations.
Moreover, the network's structure can be leveraged to construct local gradient approximations at a lower computation cost.
For instance, in the case of ODE-based network architectures, the parallel-in-time strategies~\cite{gunther2020layer} could be explored.
We do not investigate such code optimization techniques in the present work, as their effectiveness depends heavily on the particular implementation and the choice of network architecture.
However, we plan to pursue this avenue of research in the future.
Taking all aforementioned points into account, we note that the estimate $\# g_e$ considered in Table~\ref{tab:computational_cost} is fairly pessimistic.

\oldtext{The summary of computational cost and memory requirements for all aforementioned optimizers can be found in }\cancel{\Cref{tab:computational_cost}.}

\begin{table}[t]
\centering
	{
	\caption{Estimates for the number of the loss function ($\# \pazocal{L}_{\text{e}}$) and gradient evaluations ($\# g_{\text{e}}$), update cost (UC) and memory requirements (MC) per iteration and GPU device for various optimizers.}
	\label{tab:computational_cost}
	\begin{tabular}{|l|l|l|l|l|l|}
		\hline
		Method                 & $\# \pazocal{L}_{\text{e}}$                                & $\# g_{\text{e}}$                                & UC                                       & MC                                \\  \hline \hline
		Adam                   & $1$                                                        & $1$                                              & $5n$                                     & $4n$                              \\  \cline{1-5}
		L-BFGS                 & $1 + {ls}_{\text{its}}$                                    & $1$                                              & $n+ 4mn$                                 & $n+2mn$                           \\  \hline \multirow{2}{*}{ASPQN} & $1 + {ls}_{\text{its}} + $                                 & $2 + $                                           & $2n+4mn+ $                               & $2n+2mn+$                         \\
				       & $  k_{s} (\# \pazocal{L}_{\text{e}_{s}})$                  & $ k_{s} (\# g_{\text{e}_{s}})$                   & $(\sfrac{n_{s}}{n})k_{s} \text{UC}_{s}$  & $(\sfrac{n_{s}}{n}) \text{MC}_{s}$ \\ \hline
		\multirow{2}{*}{MSPQN} & $1 + {ls}_{\text{its}} + $                                 & $2 + $                                           & $2n+4mn+$                                & $2n+2mn+$                         \\
				       & $ \sum_{s=1}^{N_{sd}} k_s (\# \pazocal{L}_{\text{e}_{s}})$ & $ \sum_{s=1}^{N_{sd}} k_s (\# g_{\text{e}_{s}})$ & $k_{s} \text{UC}_{s}$                    & $(\sfrac{n_{s}}{n}) \text{MC}_{s}$ \\ \hline
	\end{tabular}
	}
\end{table}
 \section{Numerical results}
\label{sec:results}
In this section, we investigate the numerical performance of the proposed SPQN methods.
During our investigation, we monitor the loss functional~$\pazocal{L}$ as well as the quality of PINN's solution~$(u_{\texttt{NN}})$, assessed by means of the
relative $L^2$-error, given as
\begin{equation}
	\pazocal{E}_{\text{rel}}(u_{\texttt{NN}}, u^*) = \frac{ \| u_{\texttt{NN}} - u^* \|_{L^2(\Omega)} }{ \| u_{\texttt{NN}} \|_{L^2(\Omega)}}.
\end{equation}
Here, $u^*$ denotes an exact solution obtained using an analytical expression (Klein-Gordon) or solving the underlying PDE using a high-fidelity finite element solver (Burgers', Diffusion-advection, and Allen-Cahn).

\subsection{Sensitivity with respect to a varying number of subnetworks and local iterations}
\label{sec:experiments_decomp_study}
We start our numerical investigation by examining the performance of the SPQN methods with respect to a varying number of local iterations ($k_s$) and subnetworks ($N_{sd}$).
During these experiments, we consider $k_s \in \{ 10, 50, 100 \}$ and three different types of network decomposition.
In particular, we employ a decomposition into two subnetworks (the minimal decomposition), decomposition by grouping two neighboring layers into one subnetwork, and decomposition of one layer per subnetwork (the maximal layer-wise decomposition).

\begin{figure}
	\centering
	\includegraphics{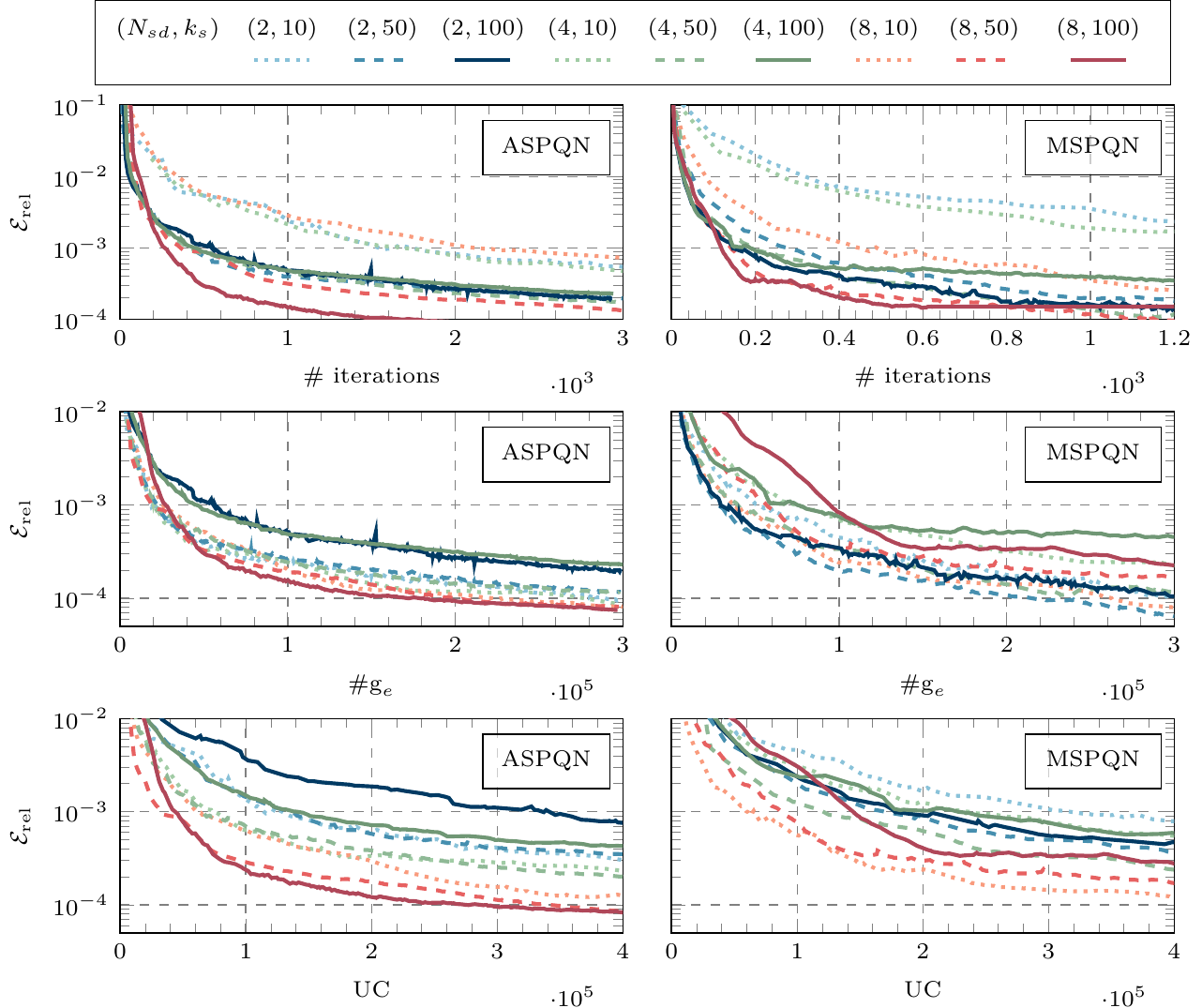}

\caption{\newtext{Burgers' equation: The mean computational cost of the ASPQN method (left column) and MSPQN method (right column).
			The results are obtained for varying numbers of local iterations ($k_s \in \{ 10, 50, 100\}$) and varying number of subdomains ($N_{sd}\in \{ 2, 4, 8 \} $). Average is obtained over $10$ independent runs.}}
	\label{fig:burgers_plots}
\end{figure}
 \begin{figure}
	\centering
\includegraphics{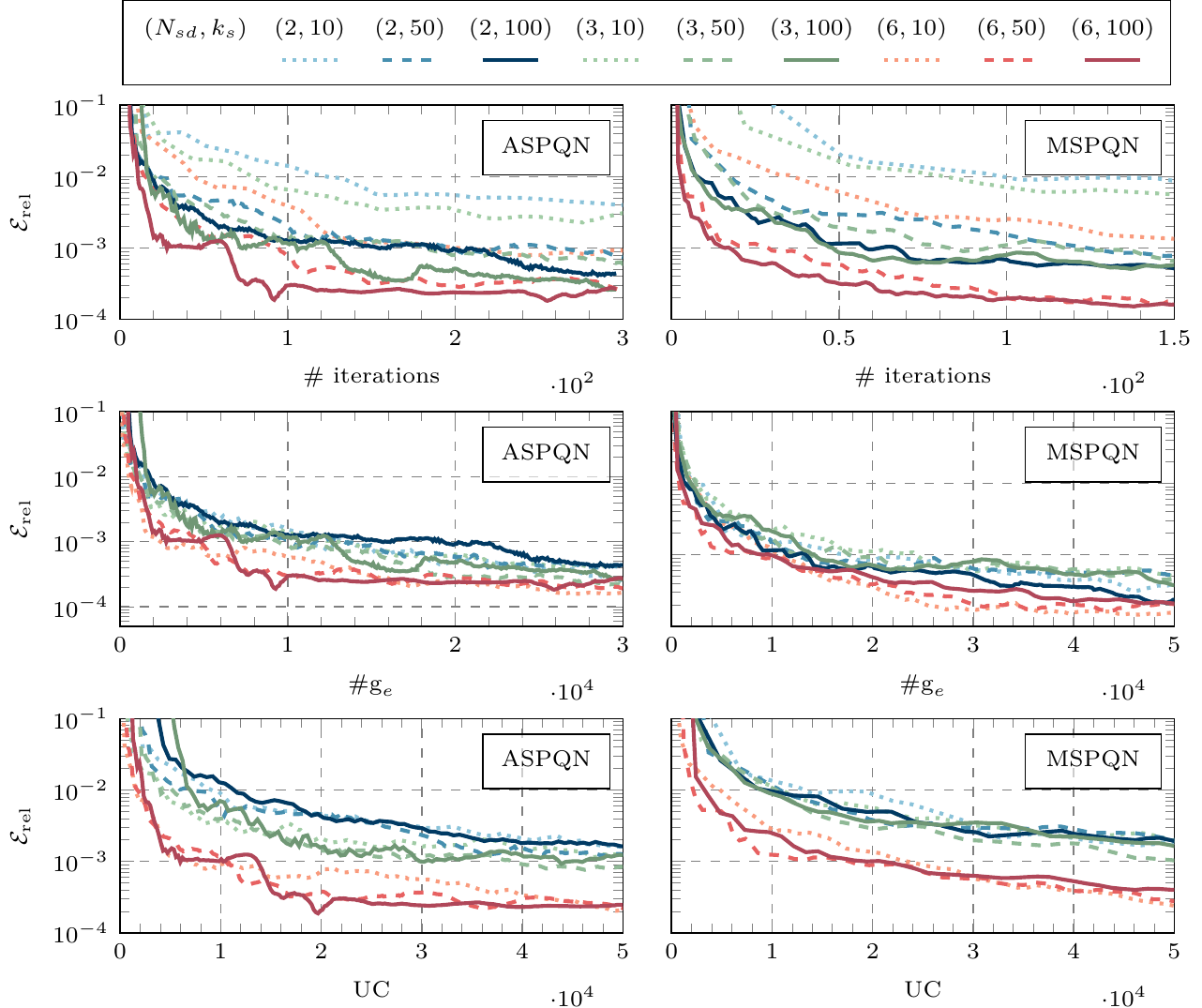}

\caption{\newtext{Klein-Gordon: The mean computational cost of the ASPQN method (left column) and MSPQN method (right column). The results are obtained for varying numbers of local iterations, i.e., ($k_s \in \{ 10, 50, 100\}$) and varying number of subdomains ($N_{sd}\in \{ 2, 3, 6 \} $). Average is obtained over $10$ independent runs.}}
	\label{fig:kleingordon_plots}
\end{figure}

\begin{figure}
	\centering
\includegraphics{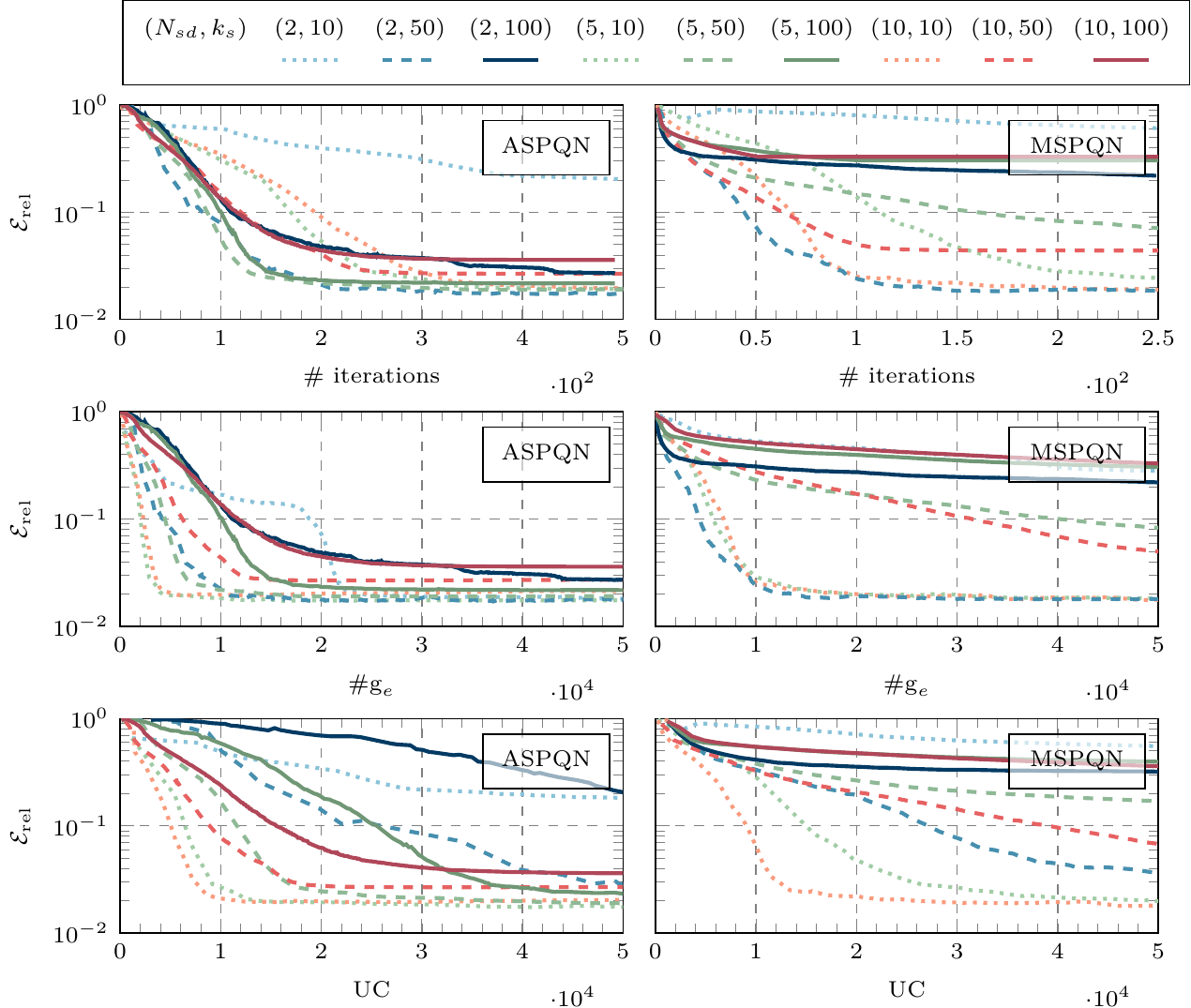}

\caption{\newtext{Advection diffusion problem: The mean computational cost of the ASPQN method (left column) and MSPQN method (right column).
			The results are obtained for varying numbers of local iterations($k_s \in \{ 10, 50, 100\}$) and varying number of subdomains ($N_{sd}\in \{ 2, 5, 10 \} $). Average is obtained over $10$ independent runs.}}
	\label{fig:transport_plots}
\end{figure}

\begin{figure}
	\centering
\includegraphics{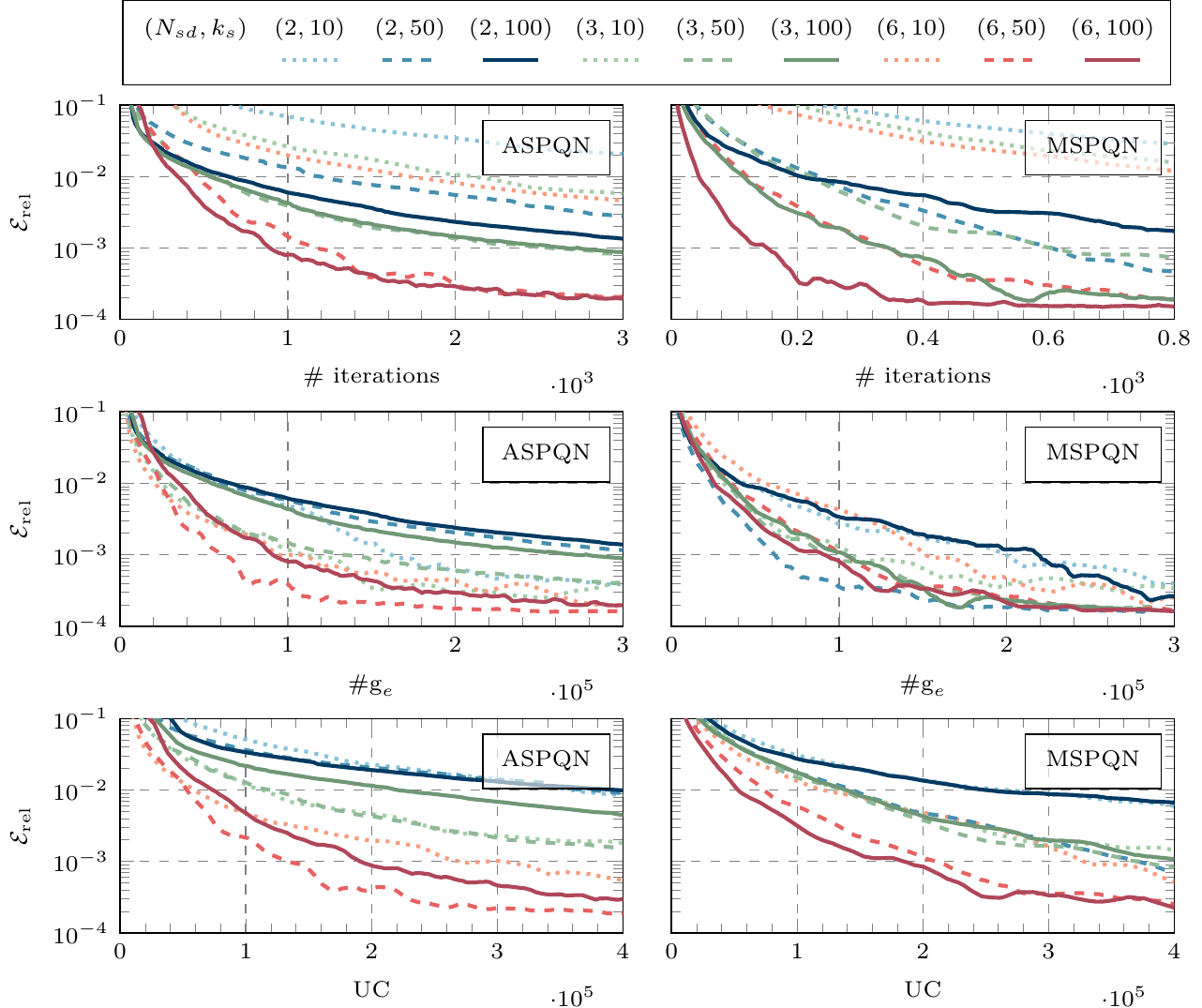}

\caption{\newtext{Allen Cahn: The mean computational cost of the ASPQN method (left column) and MSPQN method (right column).
			The results are obtained for varying numbers of local iterations ($k_s \in \{ 10, 50, 100\}$) and varying number of subdomains ($N_{sd}\in \{ 2, 3, 6 \} $). Average is obtained over $10$ independent runs.}}
	\label{fig:allencahn_plots}
\end{figure}

\cancel{\Cref{fig:aspqn_epochs}} \oldtext{illustrates the obtained results for different configurations of the ASPQN method.
More precisely, we monitor the loss~$\pazocal{L}(\boldsymbol{\theta})$ and error~$\pazocal{E}_{\text{rel}}$ as a function of number of global iterations.
As expected, we observe that the loss decreases monotonically due to the line-search globalization strategy.}
\newtext{Figures~\ref{fig:burgers_plots}-\ref{fig:allencahn_plots} demonstrate the obtained results for all four numerical examples.}
\newtext{As we} \oldtext{We also} \newtext{can} observe, \oldtext{that} the configuration of ASPQN with maximal decomposition consistently outperforms other configurations in terms of the error $\pazocal{E}_{\text{rel}}$.
The observed behavior is independent of the number of local iterations $k_s$.\oldtext{and the benchmark problem.}
This can be attributed to the fact that more extensive decoupling of the parameters allows us to obtain more localized learning rates and Hessian approximations, which in turn enables the capturing of the underlying nonlinearities more accurately.
\newtext{Moreover, we can also observe that employing the maximal layer-wise decomposition gives rise to the most efficient variants of the ASPQN in terms of $\# g_e$ and estimated UC.
The lowest $\pazocal{E}_{\text{rel}}$ is typically achieved using $50$ or $100$ local iterations.}

\cancel{\Cref{fig:aspq_geval_uc}} \oldtext{demonstrates the performance of the ASPQN method by means of the upper bound on the number of effective gradient evaluations~($\# g_e$) as well as the estimated UC, both discussed in~}\cancel{\cref{sec:implementation}.}
\oldtext{We again observe that the most efficient variants of the ASPQN method employ the maximal layer-wise decomposition.
Moreover, the lowest $\pazocal{E}_{\text{rel}}$ is typically achieved using $50$ or $100$ local iterations. }


\oldtext{In the case of the MSPQN method, the obtained results are reported in }\cancel{\cref{fig:mspqn_loss} and~\cref{fig:mspqn_cost}.}
\oldtext{As we can see from }\cancel{\Cref{fig:mspqn_loss},} \oldtext{ using a larger number of subnetworks and higher values of $k_s$ gives typically rise to a faster convergence, i.e., a fewer number of iterations is required to obtain a comparably equally accurate solution.
	By examining results reported in } \cancel{\Cref{fig:mspqn_cost},} \oldtext{ we see that the best-performing variants in terms of UC are ones with $k_s=50$ or $k_s=100$, and with the maximal layer-wise decomposition.
	This behavior is expected since the local contributions to UC decrease with an increasing number of subnetworks.
	However, the obtained results are fairly different in the context of the $\# g_e$ metric.
	As we can notice, variants with fewer subnetworks are often more effective, as the $\# g_e$ estimate does not scale proportionally with the
	size of the subnetwork, i.e., the number of trainable subnetworks' parameters. }

\newtext{In the case of the MSPQN method, the obtained results, also reported in Figures~\ref{fig:burgers_plots}-\ref{fig:allencahn_plots}, demonstrate that using a larger amount of subnetworks and higher values of $k_s$ gives typically rise to a faster convergence, i.e., a fewer number of iterations is required to obtain a comparably equally accurate solution.
The best-performing variants in terms of UC are ones with $k_s=50$ or $k_s=100$ and with the maximal layer-wise decomposition.
This behavior is expected since the local contributions to UC decrease with an increasing number of subnetworks.
In the context of the $\# g_e$ metric, the obtained results demonstrate that MSPQN variants with fewer subdomains frequently demonstrate comparable or superior efficiency in minimizing the relative error~$\pazocal{E}_{\text{rel}}$, when compared to their counterparts with full maximal layer-wise decomposition.
This is due to the fact that the $\# g_e$ estimate does not linearly scale with the subnetwork's size, i.e., the number of trainable subnetworks' parameters. }

\subsection{SPQN methods versus the state-of-the-art optimizers}
\label{sec:comparison}
\begin{figure}
	\centering
\includegraphics{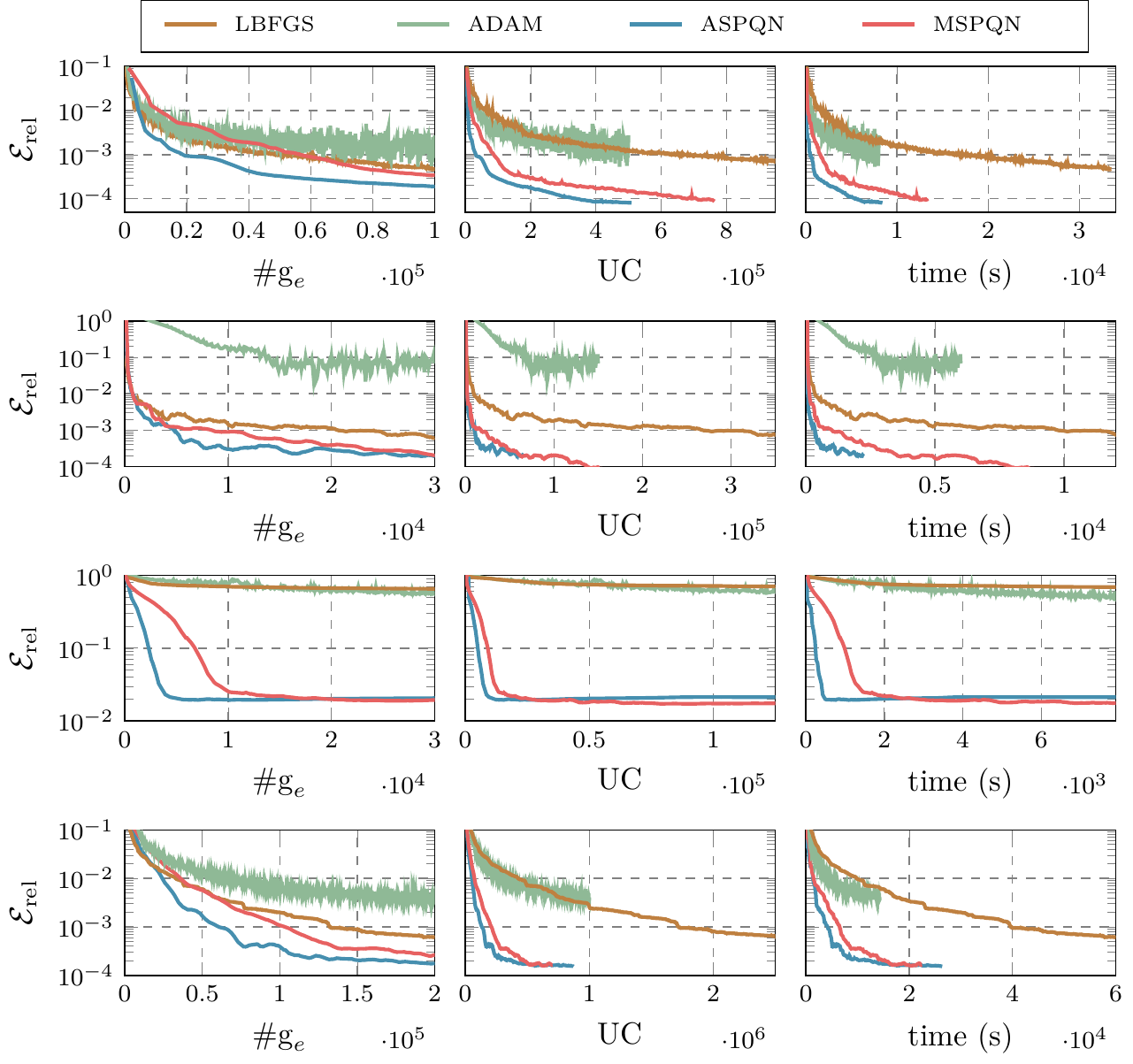}
\caption{Performance comparison of different optimizers for Burgers', Klein-Gordon, Diffusion-advection, and Allen-Cahn problems (from top to bottom). The mean relative error $\pazocal{E}_{\text{rel}}$ obtained over $10$ independent runs.}
	\label{fig:comparison}
\end{figure}
In this section, we compare the convergence and performance of the proposed SPQN methods with state-of-the-art Adam and L-BFGS optimizers.
\newtext{During these experiments, the SPQN methods consider the maximal decomposition. Moreover, we set the number of the local iterations to~ $50$. We select this particular setup as it works the best on average.
	The only exception is the advection-diffusion problem, for which the SPQN methods performed best with the small number of local iterations, in particular $k_s=10$.}
More precisely, we compare all considered optimizers by means of the $\#$$g_e$ and UC estimates.
	Moreover, we also report the execution time, although our implementation of the SPQN algorithms is not yet optimized, unlike {\sf PyTorch} implementation of Adam and L-BFGS\footnote{We have adjusted {\sf PyTorch} implementation of the L-BFGS method to take into account the momentum.}.
	\oldtext{Throughout these experiments, we configure the SPQN methods according to the specifications summarized in }\cancel{\Cref{tab:net_architectures}.}
	\oldtext{The learning rate for the Adam optimizer is chosen by carrying out a thorough hyper-parameter search.
		In particular, we choose a learning rate (reported in }\cancel{\Cref{tab:net_architectures}}\oldtext{), which allows us to obtain the most accurate PDE solution, thus a solution with the lowest~$\pazocal{E}_{\text{rel}}$.
	For the L-BFGS optimizer, we employ the same memory size ($m=3$), line-search strategy, and momentum setup as for the L-BFGS methods employed with the SPQN framework.}

	As we can see from the obtained results, reported in \Cref{fig:comparison}, the SPQN methods give rise to more accurate solutions, i.e., solutions with smaller $\pazocal{E}_{\text{rel}}$.
	The differences in the quality of an obtained solution are, on average, of one order of magnitude.
	The diffusion-advection problem is particularly interesting as the L-BFGS method stagnates for this example.
	In contrast, the variants of the SPQN optimizers are able to reach the $\pazocal{E}_{\text{rel}}$ close to $10^{-2}$.
	Numerical evidence suggests that the L-BFGS stagnation is caused by the line-search method providing very small step sizes.
	In contrast, when using the SPQN methods, a very small step size was observed only for one particular subnetwork at the beginning of the training process.
	This suggests that using the DD approach enables the usage of localized learning rates and rebalances the global nonlinearities.

	In terms of the computational cost, we can also observe from \Cref{fig:comparison}, that the SPQN methods are significantly more efficient in terms of~$\# g_e$ and UC.
	The difference is more prevalent in the case of the UC metric, as this estimate scales with the number of subnetworks.
	We can also see that the UC metric translates closely to the computational time (the last column of \Cref{fig:comparison}).
	Thus, in order to achieve a solution of comparable quality, the SPQN methods are significantly faster than Adam and L-BFGS optimizers.
	To highlight the difference, we also report the minimum mean $\pazocal{E}_{\text{rel}}$ reached by the L-BFGS optimizer and the average time required to reach that solution by all considered optimizers\footnote{Comparison with the Adam optimizer is skipped, as it never reaches $\pazocal{E}_{\text{rel}}$ as low as the L-BFGS optimizer.}
	\footnote{Comparison for the Diffusion-advection problem is skipped, as the L-BFGS optimizer stagnates for this problem.}.
	As we can see from \Cref{tab:min_error}, the computational time required by ASPQN and MSPQN is significantly lower than the computational time required by the L-BFGS method.
	In particular, the MSPQN optimizer enjoys an average speedup of a factor of $10$ across all benchmark problems.
	For the ASPQN optimizer, the speedup is even larger by an average factor of approximately $28$.
This speedup is, however, obtained with an increased amount of computational resources as the optimizer runs in parallel, i.e., using multiple GPU devices.
Thus, for a fixed amount of computational resources, MSPQN can be considered to be a more efficient optimizer.
However, the ASPQN method allows for a larger overall speedup by leveraging model parallelism.

\begin{table}[t]
	\centering
	\caption{Minimum mean $\pazocal{E}_{\text{rel}}$ reached by L-BFGS optimizer and time to solution required to reach that error for all considered optimizers.
	The network architectures and optimizer setup are the same as reported in \Cref{tab:net_architectures}.
	All optimizers are executed on a single GPU, except ASPQN, which utilizes multiple GPU devices.}
	\label{tab:min_error}
\small
	\begin{tabular}{| r | r | r | r | r  l | r |}
		\hline
		\multirow{2}{*}{{Example} } & \multirow{2}{*}{{$\pazocal{E}_{\text{rel}}$ (L-BFGS) }} & \multicolumn{5}{c|}{{Time to solution (mins)}}                                       \\  \cline{3-7}
		                            &                                                         & L-BFGS                                         & Adam & ASPQN & ($\#$ GPUs) & MSPQN  \\ \hline \hline
		Burgers'                    & $4.6 \times 10^{-4}$                                    & \ \ 558.5                                      & --   & 14.4  & (8)         & \ 40.7 \\  \hline
Klein-Gordon                & $6.1 \times 10^{-4}$                                    & \ \ 236.5                                      & --   & \ 6.8 & (6)         & \ 26.9 \\  \hline
		Allen-Cahn                  & $6.0 \times 10^{-4}$                                    & 1,001.6                                        & --   & 79.2  & (6)         & 117.5  \\  \hline
	\end{tabular}
\end{table}
 \section{Conclusion}
\label{sec:conclusion}
In this work, we proposed novel additive and multiplicative preconditioning strategies for the L-BFGS optimizer with a particular focus on PINNs applications.
The preconditioners were constructed by leveraging the layer-wise decomposition of the network, thus tailoring them to the characteristics of DNNs.
The multiplicative preconditioner was designed to enhance the convergence of the L-BFGS optimizer in the context of single GPU computing environments.
In contrast, the additive preconditioner provided a novel approach to model parallelism with a focus on multi-GPU computing architectures.
To evaluate the effectiveness of the proposed training methods, we conducted numerical experiments using benchmark problems from the field of PINNs.
A comparison with Adam and L-BFGS optimizers was performed and demonstrated that the SPQN methods can provide a substantial reduction in terms of training time.
\oldtext{In particular, using the MSPQN optimizer, we obtained an average speedup of a factor of $10$, while for the ASPQN optimizer, the average speedup factors were approximately $28$. }
\newtext{In particular, compared to L-BFGS, we obtained an approximate average speedup of $10$ for the MSPQN method and $28$ for the ASPQN method.}
Moreover, we have also demonstrated that by using the proposed SPQN methods, one can obtain a significantly more accurate solution for the underlying PDEs.

We foresee the extension of the presented work in several ways.
For instance, it would be interesting to investigate different network decompositions, as well as to incorporate overlap and coarse-level acceleration into the proposed algorithmic framework.
Furthermore, developing strategies for load balancing and optimization of the subnetwork's loss/gradient evaluations could further enhance the performance of the proposed algorithms.
\newtext{In the context of PINNs, the proposed method could be combined with the methods that consider the decomposition of the computational domain (data space).}
Moreover, the applicability of the proposed SPQN methods could be extended beyond the PINNs to a variety of different deep-learning tasks, e.g., operator learning using DeepONets~\cite{lu2021learning, lu2022comprehensive}, or natural language processing using transformers~\cite{vaswani2017attention}.
\newtextreviewerone{To this aim, the local and global optimizers would have to be adjusted to account for the noise arising from the subsampling.}
\oldtext{Lastly, in the context of PINNs, the current proposal could be combined with the methods that consider the decomposition of the computational domain (data space).}

\bibliographystyle{siamplain}
\bibliography{extracted.bib}

\begin{thebibliography}{10}

\bibitem{ben2019demystifying}
{\sc T.~Ben-Nun and T.~Hoefler}, {\em Demystifying parallel and distributed
  deep learning: An in-depth concurrency analysis}, ACM Computing Surveys
  (CSUR), 52 (2019), pp.~1--43, \url{https://doi.org/10.1145/3320060}.

\bibitem{botev2017practical}
{\sc A.~Botev, H.~Ritter, and D.~Barber}, {\em Practical {G}auss-{N}ewton
  optimisation for deep learning}, in International Conference on Machine
  Learning, PMLR, 2017, pp.~557--565, \url{https://arxiv.org/abs/1706.03662}.

\bibitem{brune2015composing}
{\sc P.~R. Brune, M.~G. Knepley, B.~F. Smith, and X.~Tu}, {\em Composing
  scalable nonlinear algebraic solvers}, SIAM Review, 57 (2015), pp.~535--565,
  \url{https://doi.org/10.1137/130936725},
  \url{https://arxiv.org/abs/1607.04254}.

\bibitem{byrd1994representations}
{\sc R.~H. Byrd, J.~Nocedal, and R.~B. Schnabel}, {\em Representations of
  quasi-{N}ewton matrices and their use in limited memory methods},
  Mathematical Programming, 63 (1994), pp.~129--156,
  \url{https://doi.org/10.1007/bf01582063}.

\bibitem{cai2002nonlinearly}
{\sc X.-C. Cai and D.~E. Keyes}, {\em {Nonlinearly preconditioned inexact
  {N}ewton algorithms}}, SIAM Journal on Scientific Computing, 24 (2002),
  pp.~183--200, \url{https://doi.org/10.1137/s106482750037620x}.

\bibitem{callens2020using}
{\sc A.~Callens, D.~Morichon, S.~Abadie, M.~Delpey, and B.~Liquet}, {\em Using
  random forest and gradient boosting trees to improve wave forecast at a
  specific location}, Applied Ocean Research, 104 (2020), p.~102339,
  \url{https://doi.org/10.1016/j.apor.2020.102339}.

\bibitem{cyr2020robust}
{\sc E.~C. Cyr, M.~A. Gulian, R.~G. Patel, M.~Perego, and N.~A. Trask}, {\em
  Robust training and initialization of deep neural networks: {A}n adaptive
  basis viewpoint}, in Mathematical and Scientific Machine Learning, PMLR,
  2020, pp.~512--536, \url{https://arxiv.org/abs/1912.04862}.

\bibitem{dean2012large}
{\sc J.~Dean, G.~Corrado, R.~Monga, K.~Chen, M.~Devin, M.~Mao, M.~Ranzato,
  A.~Senior, P.~Tucker, K.~Yang, et~al.}, {\em Large scale distributed deep
  networks}, in Advances in neural information processing systems, 2012,
  pp.~1223--1231.

\bibitem{dennis1996numerical}
{\sc J.~Dennis and R.~Schnabel}, {\em Numerical Methods for Unconstrained
  Optimization and Nonlinear Equations}, SIAM, 1996,
  \url{https://doi.org/10.1137/1.9781611971200}.

\bibitem{deuflhard2011newton}
{\sc P.~Deuflhard}, {\em {Newton} methods for nonlinear problems: affine
  invariance and adaptive algorithms}, vol.~35, Springer Science \& Business
  Media, 2011, \url{https://doi.org/10.1007/978-3-642-23899-4}.

\bibitem{dolean2022finite}
{\sc V.~Dolean, A.~Heinlein, S.~Mishra, and B.~Moseley}, {\em Finite basis
  physics-informed neural networks as a {Schwarz} domain decomposition method},
  arXiv preprint arXiv:2211.05560,  (2022).

\bibitem{dolean2023multilevel}
{\sc V.~Dolean, A.~Heinlein, S.~Mishra, and B.~Moseley}, {\em Multilevel domain
  decomposition-based architectures for physics-informed neural networks},
  arXiv preprint arXiv:2306.05486,  (2023),
  \url{https://arxiv.org/abs/2306.05486}.

\bibitem{dwivedi2021distributed}
{\sc V.~Dwivedi, N.~Parashar, and B.~Srinivasan}, {\em Distributed learning
  machines for solving forward and inverse problems in partial differential
  equations}, Neurocomputing, 420 (2021), pp.~299--316,
  \url{https://doi.org/10.1016/j.neucom.2020.09.006}.

\bibitem{weinan2021dawning}
{\sc W.~E}, {\em The dawning of a new era in applied mathematics}, Notices of
  the American Mathematical Society, 68 (2021), pp.~565--571,
  \url{https://doi.org/10.1090/noti2259}.

\bibitem{yu2018deep}
{\sc W.~E and B.~Yu}, {\em The deep {Ritz} method: a deep learning-based
  numerical algorithm for solving variational problems}, Communications in
  Mathematics and Statistics, 6 (2018), pp.~1--12,
  \url{https://doi.org/10.1007/s40304-018-0127-z},
  \url{https://arxiv.org/abs/1710.00211}.

\bibitem{fang2023ensemble}
{\sc Z.~Fang, S.~Wang, and P.~Perdikaris}, {\em Ensemble learning for physics
  informed neural networks: a gradient boosting approach}, arXiv preprint
  arXiv:2302.13143,  (2023), \url{https://arxiv.org/abs/2302.13143}.

\bibitem{Kopanicakova_2020c}
{\sc L.~Gaedke-Merzh{\"a}user, A.~l. {Kopani{\v{c}}{\'a}kov{\'a}}, and
  R.~Krause}, {\em Multilevel minimization for deep residual networks}, in
  Proceedings of French-German-Swiss Optimization Conference (FGS'2019), 2021,
  \url{https://doi.org/10.1051/proc/202171131}.

\bibitem{gavin2019levenberg}
{\sc H.~P. Gavin}, {\em The {L}evenberg-{M}arquardt algorithm for nonlinear
  least squares curve-fitting problems}, Department of Civil and Environmental
  Engineering, Duke University,  (2019), pp.~1--19.

\bibitem{gilbert1993automatic}
{\sc J.~C. Gilbert and J.~Nocedal}, {\em Automatic differentiation and the step
  computation in the limited memory {BFGS} method}, Applied mathematics
  letters, 6 (1993), pp.~47--50,
  \url{https://doi.org/10.1016/0893-9659(93)90032-i}.

\bibitem{glorot2010understanding}
{\sc X.~Glorot and Y.~Bengio}, {\em Understanding the difficulty of training
  deep feedforward neural networks}, in Proceedings of the thirteenth
  international conference on artificial intelligence and statistics, JMLR
  Workshop and Conference Proceedings, 2010, pp.~249--256.

\bibitem{golub2003separable}
{\sc G.~Golub and V.~Pereyra}, {\em Separable nonlinear least squares: the
  variable projection method and its applications}, Inverse problems, 19
  (2003), p.~R1, \url{https://doi.org/10.1088/0266-5611/19/2/201}.

\bibitem{gratton2023multilevel}
{\sc S.~Gratton, A.~Kopani\v{c}{\'a}kov{\'a}, and P.~L. Toint}, {\em Multilevel
  objective-function-free optimization with an application to neural networks
  training}, SIAM Journal on Optimization, 33 (2023), pp.~2772--2800,
  \url{https://doi.org/10.1137/23M1553455}.

\bibitem{gratton2023block}
{\sc S.~Gratton, V.~Mercier, E.~Riccietti, and P.~L. Toint}, {\em A
  block-coordinate approach of multi-level optimization with an application to
  physics-informed neural networks}, arXiv preprint arXiv:2305.14477,  (2023),
  \url{https://arxiv.org/abs/2305.14477}.

\bibitem{gu2022decomposition}
{\sc L.~Gu, W.~Zhang, J.~Liu, and X.-C. Cai}, {\em Decomposition and
  composition of deep convolutional neural networks and training acceleration
  via sub-network transfer learning}, ETNA - Electronic Transactions on
  Numerical Analysis,  (2022), pp.~157--186,
  \url{https://doi.org/10.1553/etna_vol56s157}.

\bibitem{gu2023decomposition}
{\sc L.~Gu, W.~Zhang, J.~Liu, and X.-C. Cai}, {\em Decomposition and
  preconditioning of deep convolutional neural networks for training
  acceleration}, in Domain Decomposition Methods in Science and Engineering
  XXVI, Springer, 2023, pp.~153--160.

\bibitem{gunther2020layer}
{\sc S.~G{\"u}nther, L.~Ruthotto, J.~B. Schroder, E.~C. Cyr, and N.~R. Gauger},
  {\em Layer-parallel training of deep residual neural networks}, SIAM Journal
  on Mathematics of Data Science, 2 (2020), pp.~1--23,
  \url{https://doi.org/10.1137/19m1247620}.

\bibitem{guo2015numerical}
{\sc P.~Guo, K.~Liew, and P.~Zhu}, {\em {Numerical solution of nonlinear
  Klein--Gordon equation using the element-free kp-{Ritz} method}}, Applied
  Mathematical Modelling, 39 (2015), pp.~2917--2928,
  \url{https://doi.org/10.1016/j.apm.2014.11.025}.

\bibitem{hu2022augmented}
{\sc Z.~Hu, A.~D. Jagtap, G.~E. Karniadakis, and K.~Kawaguchi}, {\em Augmented
  physics-informed neural networks ({APINNs}): A gating network-based soft
  domain decomposition methodology}, arXiv preprint arXiv:2211.08939,  (2022),
  \url{https://doi.org/10.1016/j.engappai.2023.107183},
  \url{https://arxiv.org/abs/2211.08939}.

\bibitem{hu2021extended}
{\sc Z.~Hu, A.~D. Jagtap, G.~E. Karniadakis, and K.~Kawaguchi}, {\em When do
  extended physics-informed neural networks ({XPINNs}) improve
  generalization?}, SIAM Journal on Scientific Computing, 44 (2022),
  pp.~A3158--A3182, \url{https://doi.org/10.1137/21m1447039},
  \url{https://arxiv.org/abs/2109.09444}.

\bibitem{jagtap2021extended}
{\sc A.~D. Jagtap and G.~E. Karniadakis}, {\em Extended physics-informed neural
  networks ({XPINNs}): A generalized space-time domain decomposition based deep
  learning framework for nonlinear partial differential equations}, in AAAI
  Spring Symposium: MLPS, 2021, pp.~2002--2041.

\bibitem{jagtap2020locally}
{\sc A.~D. Jagtap, K.~Kawaguchi, and G.~Em~Karniadakis}, {\em Locally adaptive
  activation functions with slope recovery for deep and physics-informed neural
  networks}, Proceedings of the Royal Society A, 476 (2020), p.~20200334,
  \url{https://doi.org/10.1098/rspa.2020.0334}.

\bibitem{jagtap2020adaptive}
{\sc A.~D. Jagtap, K.~Kawaguchi, and G.~E. Karniadakis}, {\em Adaptive
  activation functions accelerate convergence in deep and physics-informed
  neural networks}, Journal of Computational Physics, 404 (2020), p.~109136,
  \url{https://doi.org/10.1016/j.jcp.2019.109136},
  \url{https://arxiv.org/abs/1906.01170}.

\bibitem{jagtap2020conservative}
{\sc A.~D. Jagtap, E.~Kharazmi, and G.~E. Karniadakis}, {\em Conservative
  physics-informed neural networks on discrete domains for conservation laws:
  Applications to forward and inverse problems}, Computer Methods in Applied
  Mechanics and Engineering, 365 (2020), p.~113028,
  \url{https://doi.org/10.1016/j.cma.2020.113028}.

\bibitem{karniadakis2021physics}
{\sc G.~E. Karniadakis, I.~G. Kevrekidis, L.~Lu, P.~Perdikaris, S.~Wang, and
  L.~Yang}, {\em Physics-informed machine learning}, Nature Reviews Physics, 3
  (2021), pp.~422--440, \url{https://doi.org/10.1038/s42254-021-00314-5}.

\bibitem{kingma2014adam}
{\sc D.~P. Kingma and J.~Ba}, {\em Adam: A method for stochastic optimization},
  arXiv preprint arXiv:1412.6980,  (2014).

\bibitem{klawonn2023adaptive}
{\sc A.~Klawonn, M.~Lanser, and M.~Uran}, {\em Adaptive nonlinear elimination
  in nonlinear {FETI-DP} methods}, in Domain Decomposition Methods in Science
  and Engineering XXVI, Springer, 2023, pp.~337--345.

\bibitem{klawonn2023domain}
{\sc A.~Klawonn, M.~Lanser, and J.~Weber}, {\em A domain decomposition-based
  {CNN-DNN} architecture for model parallel training applied to image
  recognition problems}, arXiv preprint arXiv:2302.06564,  (2023),
  \url{https://arxiv.org/abs/2302.06564}.

\bibitem{disttraingit}
{\sc A.~Kopani{\v c}{\'a}kov{\'a} and H.~Kothari}, {\em {D}ist{T}rai{N}{N}:
  Model parallel framework for distributed training of deep neural networks.
  {G}it repository}, 2022,
  \url{{https://bitbucket.org/alena\_kopanicakova/DistTraiNN}}.

\bibitem{kopanivcakova2023nonlinear}
{\sc A.~Kopani{\v{c}}{\'a}kov{\'a}, H.~Kothari, and R.~Krause}, {\em Nonlinear
  field-split preconditioners for solving monolithic phase-field models of
  brittle fracture}, Computer Methods in Applied Mechanics and Engineering, 403
  (2023), p.~115733, \url{https://doi.org/10.1016/j.cma.2022.115733},
  \url{https://arxiv.org/abs/2203.13738}.

\bibitem{kopanivcakova2019recursive}
{\sc A.~Kopani{\v{c}}{\'a}kov{\'a}, R.~Krause, and R.~Tamstorf}, {\em
  Subdivision-based nonlinear multiscale cloth simulation}, SIAM Journal on
  Scientific Computing, 41 (2019), pp.~S433--S461,
  \url{https://doi.org/10.1137/18m1194870}.

\bibitem{kopanicakova2022globally}
{\sc A.~{Kopani{\v{c}}{\'a}kov{\'a}} and R.~Krause}, {\em Globally convergent
  multilevel training of deep residual networks}, SIAM Journal on Scientific
  Computing,  (2022), pp.~S254--S280, \url{https://doi.org/10.1137/21m1434076}.

\bibitem{kothari2022nonlinear}
{\sc H.~Kothari}, {\em Nonlinear {Schwarz} preconditioning for quasi-{Newton}
  methods}, arXiv preprint arXiv:2211.14403,  (2022),
  \url{https://arxiv.org/abs/2211.14403}.

\bibitem{kothari2022nonlinear_bounds}
{\sc H.~Kothari, A.~Kopani{\v{c}}{\'a}kov{\'a}, and R.~Krause}, {\em Nonlinear
  {S}chwarz preconditioning for nonlinear optimization problems with bound
  constraints}, arXiv preprint arXiv:2211.14780,  (2022).

\bibitem{lagaris1998artificial}
{\sc I.~E. Lagaris, A.~Likas, and D.~I. Fotiadis}, {\em Artificial neural
  networks for solving ordinary and partial differential equations}, IEEE
  transactions on neural networks, 9 (1998), pp.~987--1000,
  \url{https://doi.org/10.1109/72.712178},
  \url{https://arxiv.org/abs/physics/9705023}.

\bibitem{lee2021partition}
{\sc K.~Lee, N.~A. Trask, R.~G. Patel, M.~A. Gulian, and E.~C. Cyr}, {\em
  Partition of unity networks: deep hp-approximation}, arXiv preprint
  arXiv:2101.11256,  (2021), \url{https://doi.org/10.2172/2001532}.

\bibitem{li2019d3m}
{\sc K.~Li, K.~Tang, T.~Wu, and Q.~Liao}, {\em {D3M: A deep domain
  decomposition method for partial differential equations}}, IEEE Access, 8
  (2019), pp.~5283--5294, \url{https://doi.org/10.1109/access.2019.2957200},
  \url{https://arxiv.org/abs/1909.12236}.

\bibitem{li2020deep}
{\sc W.~Li, X.~Xiang, and Y.~Xu}, {\em Deep domain decomposition method:
  Elliptic problems}, in Mathematical and Scientific Machine Learning, PMLR,
  2020, pp.~269--286, \url{https://arxiv.org/abs/2004.04884}.

\bibitem{liu1989limited}
{\sc D.~C. Liu and J.~Nocedal}, {\em On the limited memory {BFGS} method for
  large scale optimization}, Mathematical programming, 45 (1989), pp.~503--528,
  \url{https://doi.org/10.1007/bf01589116}.

\bibitem{liu2015field}
{\sc L.~Liu and D.~E. Keyes}, {\em {Field-split preconditioned inexact {N}ewton
  algorithms}}, SIAM Journal on Scientific Computing, 37 (2015),
  pp.~A1388--A1409, \url{https://doi.org/10.1137/140970379}.

\bibitem{lu2021learning}
{\sc L.~Lu, P.~Jin, G.~Pang, Z.~Zhang, and G.~E. Karniadakis}, {\em Learning
  nonlinear operators via {D}eep{O}{N}et based on the universal approximation
  theorem of operators}, Nature Machine Intelligence, 3 (2021), pp.~218--229,
  \url{https://doi.org/10.1038/s42256-021-00302-5}.

\bibitem{lu2022comprehensive}
{\sc L.~Lu, X.~Meng, S.~Cai, Z.~Mao, S.~Goswami, Z.~Zhang, and G.~E.
  Karniadakis}, {\em A comprehensive and fair comparison of two neural
  operators (with practical extensions) based on fair data}, Computer Methods
  in Applied Mechanics and Engineering, 393 (2022), p.~114778,
  \url{https://doi.org/10.1016/j.cma.2022.114778},
  \url{https://arxiv.org/abs/2111.05512}.

\bibitem{luo2023preconditioned}
{\sc L.~Luo and X.-C. Cai}, {\em Preconditioned inexact {Newton} with learning
  capability for nonlinear system of equations}, SIAM Journal on Scientific
  Computing, 45 (2023), pp.~A849--A871.

\bibitem{lydia2019adagrad}
{\sc A.~Lydia and S.~Francis}, {\em Adagrad -- an optimizer for stochastic
  gradient descent}, Int. J. Inf. Comput. Sci, 6 (2019), pp.~566--568.

\bibitem{mahboubi2019momentum}
{\sc S.~Mahboubi, S.~Indrapriyadarsini, H.~Ninomiya, and H.~Asai}, {\em
  Momentum acceleration of quasi-{Newton} training for neural networks}, in
  PRICAI 2019: Trends in Artificial Intelligence: 16th Pacific Rim
  International Conference on Artificial Intelligence, Cuvu, Yanuca Island,
  Fiji, August 26--30, 2019, Proceedings, Part II 16, Springer, 2019,
  pp.~268--281, \url{https://doi.org/10.1007/978-3-030-29911-8_21}.

\bibitem{mcclenny2020self}
{\sc L.~McClenny and U.~Braga-Neto}, {\em Self-adaptive physics-informed neural
  networks using a soft attention mechanism}, arXiv preprint arXiv:2009.04544,
  (2020), \url{https://arxiv.org/abs/2009.04544}.

\bibitem{mercier2021coarse}
{\sc V.~Mercier, S.~Gratton, and P.~Boudier}, {\em A coarse space acceleration
  of deep-{DDM}}, arXiv preprint arXiv:2112.03732,  (2021),
  \url{https://arxiv.org/abs/2112.03732}.

\bibitem{mishra2022estimates}
{\sc S.~Mishra and R.~Molinaro}, {\em Estimates on the generalization error of
  physics-informed neural networks for approximating a class of inverse
  problems for {PDEs}}, IMA Journal of Numerical Analysis, 42 (2022),
  pp.~981--1022.

\bibitem{mnih2016asynchronous}
{\sc V.~Mnih, A.~P. Badia, M.~Mirza, A.~Graves, T.~Lillicrap, T.~Harley,
  D.~Silver, and K.~Kavukcuoglu}, {\em Asynchronous methods for deep
  reinforcement learning}, in International conference on machine learning,
  2016, pp.~1928--1937, \url{https://arxiv.org/abs/1602.01783}.

\bibitem{more1978levenberg}
{\sc J.~J. Mor{\'e}}, {\em The {L}evenberg-{M}arquardt algorithm:
  implementation and theory}, in Numerical analysis, Springer, 1978,
  pp.~105--116, \url{https://doi.org/10.1007/bfb0067700}.

\bibitem{moseley2021finite}
{\sc B.~Moseley, A.~Markham, and T.~Nissen-Meyer}, {\em Finite basis
  physics-informed neural networks ({FBPINNs}): a scalable domain decomposition
  approach for solving differential equations}, arXiv preprint
  arXiv:2107.07871,  (2021), \url{https://doi.org/10.1007/s10444-023-10065-9}.

\bibitem{newman2020train}
{\sc E.~Newman, L.~Ruthotto, J.~Hart, and B.~van Bloemen~Waanders}, {\em Train
  like a (var)pro: Efficient training of neural networks with variable
  projection}, SIAM Journal on Mathematics of Data Science, 3 (2021),
  pp.~1041--1066, \url{https://doi.org/10.1137/20m1359511},
  \url{https://arxiv.org/abs/2007.13171}.

\bibitem{o2013variable}
{\sc D.~P. O~Leary and B.~W. Rust}, {\em Variable projection for nonlinear
  least squares problems}, Computational Optimization and Applications, 54
  (2013), pp.~579--593.

\bibitem{paszke2019pytorch}
{\sc A.~Paszke, S.~Gross, F.~Massa, A.~Lerer, J.~Bradbury, G.~Chanan,
  T.~Killeen, Z.~Lin, N.~Gimelshein, L.~Antiga, et~al.}, {\em {PyTorch}: An
  imperative style, high-performance deep learning library}, Advances in neural
  information processing systems, 32 (2019).

\bibitem{rafati2018improving}
{\sc J.~Rafati and R.~F. Marcia}, {\em Improving {L-BFGS} initialization for
  trust-region methods in deep learning}, in 2018 17th IEEE International
  Conference on Machine Learning and Applications (ICMLA), IEEE, 2018,
  pp.~501--508, \url{https://doi.org/10.1109/icmla.2018.00081}.

\bibitem{raissi2019physics}
{\sc M.~Raissi, P.~Perdikaris, and G.~E. Karniadakis}, {\em Physics-informed
  neural networks: A deep learning framework for solving forward and inverse
  problems involving nonlinear partial differential equations}, Journal of
  Computational physics, 378 (2019), pp.~686--707,
  \url{https://doi.org/10.1016/j.jcp.2018.10.045}.

\bibitem{schraudolph2002fast}
{\sc N.~N. Schraudolph}, {\em Fast curvature matrix-vector products for
  second-order gradient descent}, Neural computation, 14 (2002),
  pp.~1723--1738, \url{https://doi.org/10.1162/08997660260028683}.

\bibitem{sheng2021pfnn}
{\sc H.~Sheng and C.~Yang}, {\em {PFNN}: A penalty-free neural network method
  for solving a class of second-order boundary-value problems on complex
  geometries}, Journal of Computational Physics, 428 (2021), p.~110085,
  \url{https://doi.org/10.1016/j.jcp.2020.110085},
  \url{https://arxiv.org/abs/2004.06490}.

\bibitem{sheng2022pfnn}
{\sc H.~Sheng and C.~Yang}, {\em {PFNN-2}: A domain decomposed penalty-free
  neural network method for solving partial differential equations}, arXiv
  preprint arXiv:2205.00593,  (2022),
  \url{https://doi.org/10.4208/cicp.oa-2022-0114}.

\bibitem{shukla2021parallel}
{\sc K.~Shukla, A.~D. Jagtap, and G.~E. Karniadakis}, {\em Parallel
  physics-informed neural networks via domain decomposition}, Journal of
  Computational Physics, 447 (2021), p.~110683,
  \url{https://doi.org/10.1016/j.jcp.2021.110683},
  \url{https://arxiv.org/abs/2104.10013}.

\bibitem{sirignano2018dgm}
{\sc J.~Sirignano and K.~Spiliopoulos}, {\em {DGM}: A deep learning algorithm
  for solving partial differential equations}, Journal of computational
  physics, 375 (2018), pp.~1339--1364,
  \url{https://doi.org/10.1016/j.jcp.2018.08.029},
  \url{https://arxiv.org/abs/1708.07469}.

\bibitem{stiller2020large}
{\sc P.~Stiller, F.~Bethke, M.~B{\"o}hme, R.~Pausch, S.~Torge, A.~Debus,
  J.~Vorberger, M.~Bussmann, and N.~Hoffmann}, {\em Large-scale neural solvers
  for partial differential equations}, in Driving Scientific and Engineering
  Discoveries Through the Convergence of HPC, Big Data and AI: 17th Smoky
  Mountains Computational Sciences and Engineering Conference, SMC 2020, Oak
  Ridge, TN, USA, August 26-28, 2020, Revised Selected Papers 17, Springer,
  2020, pp.~20--34, \url{https://doi.org/10.1007/978-3-030-63393-6_2}.

\bibitem{sukumar2022exact}
{\sc N.~Sukumar and A.~Srivastava}, {\em Exact imposition of boundary
  conditions with distance functions in physics-informed deep neural networks},
  Computer Methods in Applied Mechanics and Engineering, 389 (2022), p.~114333,
  \url{https://doi.org/10.1016/j.cma.2021.114333},
  \url{https://arxiv.org/abs/2104.08426}.

\bibitem{trask2022hierarchical}
{\sc N.~Trask, A.~Henriksen, C.~Martinez, and E.~Cyr}, {\em Hierarchical
  partition of unity networks: fast multilevel training}, in Mathematical and
  Scientific Machine Learning, PMLR, 2022, pp.~271--286.

\bibitem{vaswani2017attention}
{\sc A.~Vaswani, N.~Shazeer, N.~Parmar, J.~Uszkoreit, L.~Jones, A.~N. Gomez,
  L.~Kaiser, and I.~Polosukhin}, {\em {Attention is all you need}}, Advances in
  neural information processing systems, 30 (2017).

\bibitem{wang2021understanding}
{\sc S.~Wang, Y.~Teng, and P.~Perdikaris}, {\em Understanding and mitigating
  gradient flow pathologies in physics-informed neural networks}, SIAM Journal
  on Scientific Computing, 43 (2021), pp.~A3055--A3081,
  \url{https://doi.org/10.1137/20m1318043}.

\bibitem{wang2022and}
{\sc S.~Wang, X.~Yu, and P.~Perdikaris}, {\em When and why {PINNs} fail to
  train: A neural tangent kernel perspective}, Journal of Computational
  Physics, 449 (2022), p.~110768,
  \url{https://doi.org/10.1016/j.jcp.2021.110768},
  \url{https://arxiv.org/abs/2007.14527}.

\bibitem{weiss2016survey}
{\sc K.~Weiss, T.~M. Khoshgoftaar, and D.~Wang}, {\em A survey of transfer
  learning}, Journal of Big data, 3 (2016), pp.~1--40,
  \url{https://doi.org/10.1186/s40537-016-0043-6}.

\bibitem{yu2018levenberg}
{\sc H.~Yu and B.~M. Wilamowski}, {\em {L}evenberg--{M}arquardt training}, in
  Intelligent systems, CRC Press, 2018, pp.~12--1 -- 12--16,
  \url{https://doi.org/10.1201/9781315218427-12}.

\bibitem{yuan2014review}
{\sc Y.-x. Yuan}, {\em A review on subspace methods for nonlinear
  optimization}, in Proceedings of the International Congress of Mathematics,
  2014, pp.~807--827.

\end{thebibliography}

\end{document}